\documentclass[12pt, reqno]{amsart}
\usepackage{}
\usepackage{cancel}
\usepackage{stmaryrd}
\usepackage{amsmath}
\usepackage{amsfonts}
\usepackage{amssymb}
\usepackage[all,cmtip]{xy}           %xypic macro for latex2.09
\usepackage{xspace}
\usepackage{bbding}
\usepackage{txfonts}
\usepackage[shortlabels]{enumitem}
\usepackage{cite}
\usepackage{longtable}
\usepackage{makecell}
\usepackage{tikz}
\usepackage{titletoc}

\usepackage{ifpdf}
\ifpdf
 \usepackage[colorlinks,final,hyperindex]{hyperref}
\else
 \usepackage[colorlinks,final,backref=page,hyperindex,hypertex]{hyperref}
\fi
\usepackage[active]{srcltx}

\makeatletter

\begin{document}
\newcommand {\emptycomment}[1]{} %to remove paragraphs

\baselineskip=15pt
\newcommand{\nc}{\newcommand}
\newcommand{\delete}[1]{}
\nc{\mfootnote}[1]{\footnote{#1}} % Use this to show footnotes
\nc{\todo}[1]{\tred{To do:} #1}

%\delete{
\nc{\mlabel}[1]{\label{#1}}  % Use this to suppress names
\nc{\mcite}[1]{\cite{#1}}  % Use this to suppress names
\nc{\mref}[1]{\ref{#1}}  % Use this to suppress names
\nc{\meqref}[1]{\eqref{#1}} % Use this to suppress names
\nc{\mbibitem}[1]{\bibitem{#1}} % Use this to show number
%}

\delete{
\nc{\mlabel}[1]{\label{#1}  % Use the next two lines to show names
{\hfill \hspace{1cm}{\bf{{\ }\hfill(#1)}}}}
\nc{\mcite}[1]{\cite{#1}{{\bf{{\ }(#1)}}}}  % Use this lines to show names
\nc{\mref}[1]{\ref{#1}{{\bf{{\ }(#1)}}}}  % Use this lines to show names
\nc{\meqref}[1]{\eqref{#1}{{\bf{{\ }(#1)}}}} % Use this lines to show names
\nc{\mbibitem}[1]{\bibitem[\bf #1]{#1}} % Use this to show name
}

\newcommand {\comment}[1]{{\marginpar{*}\scriptsize\textbf{Comments:} #1}}
\nc{\mrm}[1]{{\rm #1}}
\nc{\id}{\mrm{id}}  \nc{\Id}{\mrm{Id}}
\nc{\admset}{\{\pm x\}\cup (-x+K^{\times}) \cup K^{\times} x^{-1}}
%%%%%%%%%%%%%%%%%%%%%%%%

\def\a{\alpha}
\def\admt{admissible to~}
\def\ad{associative D-}
\def\asi{ASI~}
\def\aybe{aYBe~}
\def\b{\beta}
\def\bd{\boxdot}
\def\btl{\blacktriangleright}
\def\btr{\blacktriangleleft}
\def\calo{\mathcal{O}}
\def\ci{\circ}
\def\d{\delta}
\def\D{\Delta}
\newcommand{\End}{\mathrm{End}}
\def\frakB{\mathfrak{B}}
\def\G{\Gamma}
\def\g{\gamma}
\def\l{\lambda}
\def\lh{\leftharpoonup}
\def\ll{\mathfrak{L}}
\def\lr{\longrightarrow}
\def\N{Nijenhuis~}
\def\o{\otimes}
\def\om{\omega}
\def\opa{\cdot_{A}}
\def\opb{\cdot_{B}}
\def\p{\psi}
\def\sadm{$S$-admissible~}
\def\r{\rho}
\def\ra{\rightarrow}
\def\rh{\rightharpoonup}
\def\rr{r^{\#}}
\def\s{\sigma}
\def\st{\star}
\def\ti{\times}
\def\tl{\triangleright}
\def\tr{\triangleleft}
\def\v{\varepsilon}
\def\vp{\varphi}
\def\vth{\vartheta}
\def\la{left Alia~}
\def\rpn{representation~}

%%%%%%%%%%%%%%%%%%%%%%%% Statements
\newtheorem{thm}{Theorem}[section]
\newtheorem{lem}[thm]{Lemma}
\newtheorem{cor}[thm]{Corollary}
\newtheorem{pro}[thm]{Proposition}
\theoremstyle{definition}
\newtheorem{defi}[thm]{Definition}
\newtheorem{ex}[thm]{Example}
\newtheorem{rmk}[thm]{Remark}
\newtheorem{pdef}[thm]{Proposition-Definition}
\newtheorem{condition}[thm]{Condition}
\newtheorem{question}[thm]{Question}
\renewcommand{\labelenumi}{{\rm(\alph{enumi})}}
\renewcommand{\theenumi}{\alph{enumi}}

\nc{\ts}[1]{\textcolor{purple}{MTS:#1}}
\nc{\zc}[1]{\color{blue}{Zhao:#1}}
\font\cyr=wncyr10

%%%%%%%%%%%%%%%%%%%%%%%%%%%%%%%%%%%%%%%%%%%%%%%%%%%%%%%%%%%%%%%%%%

 \title{Classification of three-dimensional Nijenhuis Leibniz algebras}

 \author[Ma]{Tianshui Ma \textsuperscript{1,2*}}
 \address{1. School of Mathematics and Statistics, Henan Normal University, Xinxiang 453007, China;\ ~~2. Institute of Mathematics, Henan Academy of Sciences, Zhengzhou 450046, China}
         \email{matianshui@htu.edu.cn}

 \author[Zhao]{Chan Zhao}
 \address{Chan Zhao\\ School of Mathematics and Statistics, Henan Normal University, Xinxiang 453007, China}
        \email{zhaochan2024@stu.htu.edu.cn}

 \thanks{\textsuperscript{*}Corresponding author}

 \date{\today}

 \begin{abstract}
   There are thirteen types of three-dimensional Leibniz algebras over the real field $\mathbb{R}$ based on the classification given by S. Ayupov, B. Omirov and I. Rakhimov in [Leibniz algebras: structure and classification. {\em CRC Press}, Boca Raton, FL, 2020]. In this paper, we investigate all the Nijenhuis operators on these thirteen types of three-dimensional Leibniz algebras.
 \end{abstract}

\subjclass[2020]{
17B38,  %Yang-Baxter equations and Rota-Baxter operators
17A30.  %algebras satisfying other identities
%16T25,   %Yang-Baxter equation
%16T10.   %bialgebra
%16W99,  %rings and algebras with additional structure/none of above, but in this section
%16T05,  %Hopf algebras and their applications
%17B62.   %Lie bialgebras, Lie coalgebras
%57R56,   %Topological quantum field theories
%81R60   %noncommutative geometry
}

 \keywords{Nijenhuis operator, Leibniz algebra.}

 \maketitle

% \vspace{-1.2cm}

% \tableofcontents

%\vspace{-1cm}

 \allowdisplaybreaks

\section{Introduction and preliminaries}
 As a ``non-commutative" analogue of Lie algebras, a (left) {\bf Leibniz algebra} (\cite{B,L,LP}) is a pair $(\ll, [,])$ consisting of a vector space $\ll$ and a bilinear map $[,]: \ll\o \ll\lr \ll$ (write $[,](x\o y)=[x, y]$) such that for all $x, y, z\in \ll$,
 \begin{eqnarray*}
 &[x, [y, z]]=[[x, y], z]+[y, [x, z]].&\label{eq:l}
 \end{eqnarray*}  
 A {\bf \N Leibniz algebra} (\cite{CGM1}) is a pair $((\ll, [,]), N)$, where $(\ll, [,])$ is a Leibniz algebra and $N:\ll\lr \ll$ is a linear map such that for all $x, y\in \ll$,
 \begin{eqnarray*}
 [N(x), N(y)]+N^2([x, y])=N([N(x), y])+N([x, N(y)]).
 \end{eqnarray*}
 In this case, $N$ is called a {\bf \N operator} on $(\ll, [,])$. Leibniz algebras and \N operators have attracted the attention of many researchers. 
 
 More recently, in \cite{GD}, Guo and Das explored the concept of generalized Reynolds operators on Leibniz algebras as an extension of twisted Poisson structures, and their investigation is grounded in the Loday-Pirashvili cohomology of an induced Leibniz algebra. In \cite{MS}, Mondal and Saha discussed the relationship of Nijenhuis operators with Rota-Baxter operators and modified Rota-Baxter operators on Leibniz algebras and then considered the cohomology and deformation theory. In \cite{LMW}, Li, Ma and Wang found that Leibniz algebras are closely related to Nijenhuis operators, and in \cite{MSZ}, Ma, Sun and Zheng investigated the bialgebraic structures on Nijenhuis Leibniz algebras.
  
 In this paper, as a continuation of \cite{LMW,MSZ}, we investigate all the Nijenhuis operators on these thirteen types of three-dimensional Leibniz algebras over the real field $\mathbb{R}$, which are based on the classification given by S. Ayupov, B. Omirov and I. Rakhimov in \cite{AOR}. This is an important step towards achieving the classification of all three-dimensional Nijenhuis Leibniz bialgebras.
 \smallskip

 \section{Three dimensional Nijenhuis Leibniz algebras} \label{se:3} According to the classification given by S. Ayupov, B. Omirov and I. Rakhimov in \cite{AOR}, there are thirteen kinds of Leibniz algebras of dimension three. We now give all Nijenhuis operators on these three-dimensional Leibniz algebras. Let $(\ll, [,])$ be a Leibniz algebra of dimension three with basis $\{e, f, g\}$. Define a linear map $N: \ll\lr \ll$ by
 \begin{center}
 $\left\{
\begin{array}{l}
 N(e)=k_1 e+k_2 f+k_3 g\\
 N(f)=\ell_1 e+\ell_2 f+\ell_3 g\\
 N(g)=p_1 e+p_2 f+p_3 g\\
 \end{array}
 \right.$,
 \end{center}
 where $k_i, \ell_i, p_i, i=1, 2, 3$ are parameters. Then on the basis $\{e, f, g\}$, there are one-to-one correspondence as follows:
 \begin{center}$N\longleftrightarrow$ $\left(
 \begin{array}{ccc}
 k_1 & k_2 & k_3 \\
 \ell_1 & \ell_2 & \ell_3 \\
 p_1 & p_2 & p_3 \\
 \end{array}
 \right).$ 
 \end{center}

 \begin{thm}\label{thm:3nb1} Let $(\ll, [,])$ be a Leibniz algebra and $[,]$ is given by
 \begin{center}
 \begin{tabular}{r|rrr}
 $[,]$ & $e$  & $f$ & $g$ \\
 \hline
 $e$ & $0$  & $0$ & $0$ \\
 $f$ & $0$  & $e$ & $f$  \\
 $g$ & $-2e$  & $-f$ & $0$\\
 \end{tabular}.
 \end{center}\vskip1mm
 Then all the Nijenhuis operators on $(\ll, [,])$ are given as follows:
 \begin{center}$(1)~\left(
 \begin{array}{ccc}
 k_1 & 0 & 0 \\
 0 & k_1 & \ell_3 \\
 0 & 0 & p_3 \\
 \end{array}
 \right),$ \quad
 $(2)~\left(
 \begin{array}{ccc}
 k_1 & 0 & 0 \\
 \ell_1 & \ell_2 & \frac{(k_1-\ell_2)^2}{2 \ell_1}\\
 0 & 0 & k_1 \\
 \end{array}
 \right),~(\ell_1\neq 0),$ 
 \end{center}
 \begin{center}
 $(3)~\left(
 \begin{array}{ccc}
 k_1 & 0 & 0 \\
 0 & k_1 & 0 \\
 p_1 & p_2  & \frac{p_2^2}{ 2 p_1} - k_1 \\
 \end{array}
 \right),~(p_1\neq 0)$\quad
 $(4)\left(
 \begin{array}{ccc}
 -\frac{p_2 \ell_1}{p_1}+ \ell_2 & 0  & 0 \\
 \ell_1 & \ell_2 & \frac{p_2^2 \ell_1}{ 2 p_1^2}\\\
 p_1 & p_2 & \frac{p_2(p_2 - 2\ell_1)}{ 2 p_1} + \ell_2 \\
 \end{array}
 \right) ~(p_1\neq 0, \ell_1\neq 0).$ 
 \end{center} 
 \end{thm}

 \begin{proof} $N$ is Nijenhuis operator on $(\ll, [,])$ if and only if
 \begin{center}
 $\left\{
 \begin{array}{l}
 k_2 = k_3 = 0\\
 (k_1-\ell_2)^2-2 \ell_1 \ell_3=0\\
 - \ell_1 p_3 + \ell_2 p_2 + \ell_1 k_1- \ell_3 p_1-p_2 k_1 = 0 \\
 -2 p_1 p_3 - p_2^2 + 2p_1 k_1 = 0\\
 \end{array}
 \right.$.
 \end{center}

 \begin{enumerate}[(I)]
 \item $p_1 = 0$, then $\left\{
 \begin{array}{l}
 k_2 = k_3 = p_1 = 0\\
 (k_1-\ell_2)^2-2 \ell_1 \ell_3=0\\
 \ell_1(p_3-k_1)= 0 \\
 p_2 = 0\\
 \end{array}
 \right.$. \smallskip
 \begin{enumerate}
 \item [(IA)] $\ell_1 = 0$, then $\left\{
 \begin{array}{l}
 k_2 = k_3 = p_1 = \ell_1 = p_2 = 0\\
 \ell_2 = k_1\\
 \end{array}
 \right.$. That is to say, 
 \begin{center}$\left\{
 \begin{array}{l}
 N(e)=k_1 e\\
 N(f)=k_1 f + \ell_3 g\\
 N(g)=p_3 g\\
 \end{array}
 \right.$.
 \end{center}\smallskip
 \item [(IB)] \label{it:IB} $\ell_1 \neq 0$, we have
 $\left\{
 \begin{array}{l}
 k_2 = k_3 = p_1 = p_2 = 0\\
 \ell_1 \neq 0\\
 p_3 = k_1\\
 \ell_3=\frac{(k_1-\ell_2)^2}{2 \ell_1}\\
 \end{array}
 \right.$. That is to say,\smallskip
 \begin{center}$\left\{
 \begin{array}{l}
 N(e)=k_1 e\\
 N(f)=\ell_1 e + \ell_2 f + \frac{(k_1-\ell_2)^2}{2 \ell_1} g\\
 N(g)=k_1 g\\
 \end{array}
 \right.~(\ell_1 \neq 0)$.
 \end{center}
 \end{enumerate}\smallskip
 \item $p_1\neq 0$, then $\left\{
 \begin{array}{l}
 k_2 = k_3 = 0\\
 p_1\neq 0\\
 (k_1-\ell_2)^2-2 \ell_1 \ell_3=0\\
 - \ell_1 p_3 + \ell_2 p_2 + \ell_1 k_1- \ell_3 p_1-p_2 k_1 = 0 \\
 p_3 = \frac{p_2^2}{2 p_1} +k_1 \\
 \end{array}
 \right.$. \smallskip
 \begin{enumerate}
 \item [(IIA)] $\ell_1=0$, then
 $\left\{
 \begin{array}{l}
 k_2 = k_3 = \ell_1 = \ell_3 =0\\
 p_1\neq 0\\
 \ell_2=k_1\\
 p_3 = \frac{p_2^2}{2 p_1} +k_1 \\
 \end{array}
 \right.$. That is to say,\smallskip
 \begin{center}$\left\{
 \begin{array}{l}
 N(e)= k_1 e\\
 N(f)= k_1 f \\
 N(g)=p_1 e + p_2 f + (\frac{p_2^2}{2 p_1} +k_1) g\\
 \end{array}
 \right. ~(p_1\neq 0)$.
 \end{center}\smallskip
 \item [(IIB)] $\ell_1\neq 0$, then $\left\{
 \begin{array}{l}
 k_2 = k_3 = 0\\
 p_1\neq 0, \ell_1\neq 0\\
 k_1 = -\frac{p_2 \ell_1}{p_1}+ \ell_2\\
 \ell_3 = \frac{p_2^2 \ell_1}{ 2 p_1^2}\\\
 p_3 = \frac{p_2(p_2 - 2\ell_1)}{ 2 p_1} + \ell_2 \\
 \end{array}
 \right.$. That is to say,\smallskip
 \begin{center}$\left\{
 \begin{array}{l}
 N(e)=(-\frac{p_2 \ell_1}{p_1}+ \ell_2) e\\
 N(f)=\ell_1 e+ \ell_2 f+\frac{p_2^2 \ell_1}{ 2 p_1^2} g\\
 N(g)=p_1 e+ p_2 f+(\frac{p_2(p_2 - 2\ell_1)}{ 2 p_1} + \ell_2) g\\
 \end{array}
 \right.~(p_1 \neq 0, \ell_1 \neq 0)$.
 \end{center}\qedhere
 \end{enumerate}
 \end{enumerate}
 \end{proof}

 \begin{thm}\label{thm:3nb2} Let $(\ll, [,])$ be a Leibniz algebra and $[,]$ is given by
 \begin{center}
 \begin{tabular}{r|rrr}
 $[,]$ & $e\quad$  & $f$ & $g$ \\
 \hline
 $e$ & $0\quad$  & $0$ & $0$ \\
 $f$ & $0\quad$  & $0$ & $0$  \\
 $g$ & $e+f$  & $0$ & $e$  \\
 \end{tabular}.
 \end{center}\smallskip
 Then all the Nijenhuis operators on $(\ll, [,])$ are given as follows:
 \begin{center}$(1)~\left(
 \begin{array}{ccc}
 k_1 & 0 & 0 \\
 0 & k_1 & 0 \\
 p_1 & p_2 & k_1 \\
 \end{array}
 \right),$ \quad
 $(2)~\left(
 \begin{array}{ccc}
 k_1 & 0 & 0 \\
 0 & \ell_2 & 0\\
 \ell_2-k_1 & 0 & \ell_2 \\
 \end{array}
 \right),~(\ell_2\neq k_1),$  
 \end{center}\smallskip

 \begin{center}
 $(3)~\left(
 \begin{array}{ccc}
 2k_2+\ell_2 & k_2 & 0 \\
 -k_2 & \ell_2 & 0 \\
 p_1 & p_2  & k_2+\ell_2 \\
 \end{array}
 \right),~(k_2\neq 0)$,\quad
 $(4)\left(
 \begin{array}{ccc}
 k_1 & k_2 & 0 \\
 0 & \ell_2 & 0\\\
 k_1-\ell_2 & p_2 & \ell_2\\
 \end{array}
 \right) ~(k_2\neq 0).$ 
 \end{center} 
 \end{thm}

 \begin{proof} $N$ is Nijenhuis operator on $(\ll, [,])$ if and only if
 \begin{center}
 $\left\{
 \begin{array}{l}
 k_3 = \ell_3 = 0\\
 (p_3-k_1)^2+p_1(p_3- k_1- \ell_1)+k_2 \ell_1=0\\
 p_1(p_3- k_2- \ell_2)+k_2(k_1+ \ell_2- 2 p_3)=0\\
 \ell_1(p_3-k_1-\ell_1)=0\\
 \ell_1(p_3-k_2-\ell_2)=0\\
 (k_2+\ell_2-k_1)(\ell_2-p_3)+\ell_1 k_2=0\\
 \end{array}
 \right.$.
 \end{center}
 \begin{enumerate}[(I)]
 \item $p_3= k_2+ \ell_2$, then
 $\left\{
 \begin{array}{l}
 k_3 = \ell_3 = 0\\
 (p_3-k_1)^2+p_1(p_3- k_1- \ell_1)+k_2 \ell_1=0\\
 k_2(k_1+ \ell_2- 2 p_3)=0\\
 \ell_1(p_3-k_1-\ell_1)=0\\
 (p_3-k_1)(\ell_2-p_3)+\ell_1 k_2=0\\
 p_3= k_2+ \ell_2\\
 \end{array}
 \right.$. \smallskip

 \begin{enumerate}
 \item [(IA)] $k_2=0$, then
 $\left\{
 \begin{array}{l}
 k_3 = \ell_3 = k_2 = 0\\
 (p_3-k_1)^2+p_1(p_3- k_1- \ell_1)=0\\
 \ell_1(p_3-k_1-\ell_1)=0\\
 p_3= \ell_2 \\
 \end{array}
 \right.$. \smallskip

 \begin{enumerate}
 \item [(IA1)] $p_3= k_1+ \ell_1$, then $\left\{
 \begin{array}{l}
 k_3 = \ell_3 = k_2 = \ell_1 = 0\\
 p_3 = k_1 = \ell_2\\
 \end{array}
 \right. $. That is to say,\smallskip
 \begin{center}$\left\{
 \begin{array}{l}
 N(e)=k_1 e\\
 N(f)=k_1 f\\
 N(g)=p_1 e+p_2 f+k_1 g\\
 \end{array}
 \right. $.
 \end{center}\smallskip

 \item [(IA2)] $p_3\neq k_1+ \ell_1$, we have $\left\{
 \begin{array}{l}
 k_3 = \ell_3 = k_2 = \ell_1 = 0\\
 p_1=k_1-\ell_2\\
 p_3= \ell_2 \\
 \ell_2\neq k_1\\
 \end{array}
 \right.$. That is to say,\smallskip
 \begin{center}$\left\{
 \begin{array}{l}
 N(e)=k_1 e\\
 N(f)=\ell_2 f\\
 N(g)=(k_1-\ell_2) e+p_2 f+\ell_2 g\\
 \end{array}
 \right.~(\ell_2\neq k_1)$.
 \end{center}
 \end{enumerate}

  \item [(IB)] $k_2\neq 0$, we have
 $\left\{
 \begin{array}{l}
 k_3 = \ell_3 = 0\\
 (p_3-k_1)^2+p_1(p_3- k_1- \ell_1)+k_2 \ell_1=0\\
 k_1+ \ell_2= 2 p_3\\
 \ell_1(p_3-k_1-\ell_1)=0\\
 (p_3-k_1)(\ell_2-p_3)+\ell_1 k_2=0\\
 p_3= k_2+ \ell_2\\
 k_2\neq 0\\
 \end{array}
 \right. $.\smallskip

 then $\left\{
 \begin{array}{l}
 k_3 = \ell_3 = 0\\
 k_1=2 k_2+ \ell_2\\
 \ell_1= -k_2\\
 p_3= k_2+ \ell_2\\
 k_2\neq 0\\
 \end{array}
 \right. $. That is to say,\smallskip
 \begin{center}$\left\{
 \begin{array}{l}
 N(e)=(2 k_2+ \ell_2) e+ k_2 f\\
 N(f)=-k_2 e+k_1 f\\
 N(g)=p_1 e+p_2 f+ (k_2+ \ell_2) g\\
 \end{array}
 \right. ~(k_2\neq 0)$.
 \end{center}
 \end{enumerate}\smallskip

 \item $p_3\neq k_2+ \ell_2$, then $\left\{
 \begin{array}{l}
 k_3 = \ell_3 = \ell_1 = 0\\
 (p_3-k_1+p_1)(p_3- k_1)=0\\
 p_1(p_3- k_2- \ell_2)+k_2(k_1+ \ell_2- 2 p_3)=0\\
 (k_2+\ell_2-k_1)(\ell_2-p_3)=0\\
 p_3\neq k_2+ \ell_2\\
 \end{array}
 \right. $. \smallskip

 \begin{enumerate}
 \item [(IIA)] $p_3=\ell_2$, then
 $\left\{
 \begin{array}{l}
 k_3 = \ell_3 = \ell_1 = 0\\
 p_1= k_1- \ell_2\\
 p_3=\ell_2\\
 k_2\neq 0\\
 \end{array}
 \right.$. That is to say,\smallskip
 \begin{center}$\left\{
 \begin{array}{l}
 N(e)=k_1 e+ k_2 f\\
 N(f)=\ell_2 f\\
 N(g)=(k_1- \ell_2) e+p_2 f+\ell_2 g\\
 \end{array}
 \right.~(k_2\neq 0)$.
 \end{center}\smallskip

 \item [(IIB)] $p_3 \neq \ell_2$, then
 $\left\{
 \begin{array}{l}
 k_3 = \ell_3 = \ell_1 = 0\\
 p_3 = \ell_2\\
 k_1=k_2+\ell_2\\
 p_3 \neq \ell_2, p_3\neq k_2+ \ell_2\\
 \end{array}
 \right.$. This assumption is not valid.
 \end{enumerate}\qedhere
 \end{enumerate}
 \end{proof}

 \begin{thm}\label{thm:3nb3}  Let $(\ll, [,])$ be a Leibniz algebra and $[,]$ is given by
 \begin{center}
 \begin{tabular}{r|rrr}
 $[,]$ & $e$  & $f$ & $g$ \\
 \hline
 $e$ & $0$  & $0$ & $0$ \\
 $f$ & $0$  & $0$ & $f$  \\
 $g$ & $\a e$  & $-f$ & $0$  \\
 \end{tabular},
 where $0\neq\a\in \mathfrak{R}$ (real field).
 \end{center}\vskip1mm
 Then all the Nijenhuis operators on $(\ll, [,])$ are given as follows:
 \begin{center}$(1)~\left(
 \begin{array}{ccc}
 k_1 & 0 & 0 \\
 \ell_1 & \ell_2 & 0 \\
 p_1 & p_2 & k_1 \\
 \end{array}
 \right),$ \quad
 $(2)~\left(
 \begin{array}{ccc}
 k_1 & 0 & 0 \\
 0 & \ell_2 & \ell_3\\
 0 & p_2 & k_1 \\
 \end{array}
 \right),~(\ell_3\neq 0),$  
 \end{center}

 \begin{center}
 $(3)~\left(
 \begin{array}{ccc}
 k_1 & 0 & 0 \\
 0 & \ell_2 & \ell_3\\
 0 & p_2  & p_3 \\
 \end{array}
 \right),~(p_3\neq k_1)$\quad
 $(4)\left(
 \begin{array}{ccc}
 k_1 & k_2  & 0 \\
 0 & \ell_2 & 0\\\
 0 & p_2 & \ell_2 \\
 \end{array}
 \right) ~(k_2\neq 0).$
 \end{center} 
 \end{thm}

 \begin{proof} $N$ is Nijenhuis operator on $(\ll, [,])$ if and only if
 \begin{center}
 $\left\{
 \begin{array}{l}
 k_3 =0\\
 k_2 \ell_3=0\\
 k_2 \ell_1=0\\
 k_2 (p_3-\ell_2)=0 \\
 \ell_3 \ell_1=0\\
 \ell_3 p_1(1+\a)+\ell_1(k_1-p_3)=0\\
 \ell_3 p_1+\ell_1(k_1-p_3)(1+\a)=0\\
 k_2 p_1=0\\
 p_1 (p_3-k_1)=0 \\
 \end{array}
 \right.$.
 \end{center}
 \begin{enumerate}[(I)]
 \item $k_2=0$, then $\left\{
 \begin{array}{l}
 k_3 = k_2 = 0\\
 \ell_3 \ell_1=0\\
 \ell_3 p_1(1+\a)+\ell_1(k_1-p_3)=0\\
 \ell_3 p_1+\ell_1(k_1-p_3)(1+\a)=0\\
 p_1 (p_3-k_1)=0 \\
 \end{array}
 \right.$. \smallskip
 \begin{enumerate}
 \item [(IA)] $p_3=k_1$, then $\left\{
 \begin{array}{l}
 k_3 = k_2 = 0\\
 \ell_3 \ell_1=0\\
 \ell_3 p_1=0\\
 p_3=k_1\\
 \end{array}
 \right.$. \smallskip
 \begin{enumerate}
 \item [(IA1)] $\ell_3=0$, then $\left\{
 \begin{array}{l}
 k_3 = k_2 = \ell_3 = 0\\
 p_3=k_1\\
 \end{array}
 \right.$. That is to say,\smallskip
 \begin{center}$\left\{
 \begin{array}{l}
 N(e)=k_1 e\\
 N(f)=\ell_1 e+\ell_2 f\\
 N(g)=p_1 e+p_2 f+k_1 g\\
 \end{array}
 \right.$.
 \end{center}\smallskip
 \item [(IA2)] $\ell_3\neq 0$, we have $\left\{
 \begin{array}{l}
 k_3 = k_2 = \ell_1 = p_1 = 0\\
 p_3=k_1\\
 \ell_3\neq 0\\
 \end{array}
 \right.$. That is to say,\smallskip
 \begin{center}$\left\{
 \begin{array}{l}
 N(e)=k_1 e\\
 N(f)=\ell_2 f+\ell_3 g\\
 N(g)=p_2 f+k_1 g\\
 \end{array}
 \right.~(\ell_3\neq 0)$.
 \end{center}
 \end{enumerate}\smallskip
 \item [(IB)]  $p_3\neq k_1$, we have
 $\left\{
 \begin{array}{l}
 k_3 = k_2 = p_1 = \ell_1 = 0\\
 p_3\neq k_1\\
 \end{array}
 \right.$. That is to say,\smallskip
 \begin{center}$\left\{
 \begin{array}{l}
 N(e)=k_1 e\\
 N(f)=\ell_2 f+\ell_3 g\\
 N(g)=p_2 f+p_3 g\\
 \end{array}
 \right.~(p_3\neq k_1)$.
 \end{center}
 \end{enumerate}
 \item $k_2\neq0$, then $\left\{
 \begin{array}{l}
 k_3 = \ell_3 = \ell_1 = p_1 = 0\\
 p_3=\ell_2 \\
 k_2\neq0\\
 \end{array}
 \right.$. That is to say,\smallskip
 \begin{center}$\left\{
 \begin{array}{l}
 N(e)=k_1 e+k_2 f\\
 N(f)=\ell_2 f\\
 N(g)=p_2 f+\ell_2 g\\
 \end{array}
 \right.~(k_2\neq0)$.
 \end{center}\qedhere
 \end{enumerate}
 \end{proof}

 \begin{thm}\label{thm:3nb4}  Let $(\ll, [,])$ be a Leibniz algebra and $[,]$ is given by
 \begin{center}
 \begin{tabular}{r|rrr}
 $[,]$ & $e$  & $f$ & $g$ \\
 \hline
 $e$ & $0$  & $0$ & $0$ \\
 $f$ & $0$  & $0$ & $f$  \\
 $g$ & $0$  & $-f$ & $e$  \\
 \end{tabular}.
 \end{center} 
 Then all the Nijenhuis operators on $(\ll, [,])$ are given as follows:\smallskip
 \begin{center}$(1)~\left(
 \begin{array}{ccc}
 k_1 & 0 & 0 \\
 \ell_1 & \ell_2 & 0 \\
 p_1 & p_2 & k_1 \\
 \end{array}
 \right),$ \quad
 $(2)~\left(
 \begin{array}{ccc}
 k_1 & k_2 & 0 \\
 0 & k_1 & 0\\
 p_1 & p_2 & k_1 \\
 \end{array}
 \right),~(k_2\neq 0).$
 \end{center} 
 \end{thm}

 \begin{proof} $N$ is Nijenhuis operator on $(\ll, [,])$ if and only if
 \begin{center}
 $\left\{
 \begin{array}{l}
 k_3 = \ell_3 = 0\\
 k_2(p_3-\ell_2) =0\\
 k_2 \ell_1 =0 \\
 p_3= k_1\\
 k_2(k_1+\ell_2-2 p_3)=0\\
 \end{array}
 \right.$.
 \end{center}
 \begin{enumerate}[(I)]
 \item $k_2=0$, then $\left\{
 \begin{array}{l}
 k_3 = \ell_3 = k_2 = 0\\
 p_3= k_1\\
 \end{array}
 \right.$. That is to say,\smallskip
 \begin{center}$\left\{
 \begin{array}{l}
 N(e)=k_1 e\\
 N(f)=\ell_1 e+\ell_2 f\\
 N(g)=p_1 e+p_2 f+k_1 g\\
 \end{array}
 \right.$.
 \end{center}\smallskip
 \item $k_2\neq0$, then $\left\{
 \begin{array}{l}
 k_3 = \ell_3 = \ell_1 = 0\\
 p_3 = \ell_2 = k_1\\
 k_2\neq0\\
 \end{array}
 \right.$. That is to say,\smallskip
 \begin{center}$\left\{
 \begin{array}{l}
 N(e)=k_1 e+k_2 f\\
 N(f)=k_1 f\\
 N(g)=p_1 e+p_2 f+k_1 g\\
 \end{array}
 \right.~(k_2\neq0)$.
 \end{center}
 \end{enumerate}
 \end{proof}

 \begin{thm}\label{thm:3nb5}  Let $(\ll, [,])$ be a Leibniz algebra and $[,]$ is given by
 \begin{center}
 \begin{tabular}{r|rrr}
 $[,]$ & $e$  & $f$ & $g$ \\
 \hline
 $e$ & $0$  & $0$ & $0$ \\
 $f$ & $0$  & $e$ & $0$  \\
 $g$ & $0$  & $0$ & $e$  \\
 \end{tabular}.
 \end{center}\vskip1mm
 Then all the Nijenhuis operators on $(\ll, [,])$ are given as follows:
 \begin{center}$\left(
 \begin{array}{ccc}
 k_1 & 0 & 0 \\
 \ell_1 & k_1 & 0 \\
 p_1 & 0 & k_1 \\
 \end{array}
 \right).$
 \end{center} 
 \end{thm}

 \begin{proof} $N$ is Nijenhuis operator on $(\ll, [,])$ if and only if
 \begin{center}
 $\left\{
 \begin{array}{l}
 k_2 =k_3= 0\\
 (k_1-\ell_2)^2+\ell_3^2=0\\
 (p_3-k_1)^2+p_2^2=0 \\
 p_2(\ell_2-k_1)=0 \\
 \end{array}
 \right.$.
 \end{center} then $\left\{
 \begin{array}{l}
 k_2 = k_3 = \ell_3 = p_2 = 0\\
 k_1 = \ell_2 = p_3\\
 \end{array}
 \right.$. 
 That is to say,  
 \begin{center}$\left\{
 \begin{array}{l}
 N(e)=k_1 e\\
 N(f)=\ell_1 e+k_1 f\\
 N(g)=p_1 e+k_1 g\\
 \end{array}
 \right.$.
 \end{center} 
 \end{proof}

 \begin{thm}\label{thm:3nb6}  Let $(\ll, [,])$ be a Leibniz algebra and $[,]$ is given by
 \begin{center}
 \begin{tabular}{r|rrr}
 $[,]$ & $e$  & $f$ & $g$ \\
 \hline
 $e$ & $0$  & $0$ & $0$ \\
 $f$ & $0$  & $e$ & $0$  \\
 $g$ & $0$  & $0$ & $-e$  \\
 \end{tabular}.
 \end{center}\vskip1mm
 Then all the Nijenhuis operators on $(\ll, [,])$ are given as follows:
 \begin{center}$(1)~\left(
 \begin{array}{ccc}
 k_1 & 0 & 0 \\
 \ell_1 & \ell_2 & \ell_2-k_1 \\
 p_1 & p_3-k_1 & p_3 \\
 \end{array}
 \right),$  \quad %~(\ell_2\neq k_1),
 $(2)~\left(
 \begin{array}{ccc}
 k_1 & 0 & 0 \\
 \ell_1 & k_1 & 0 \\
 p_1 & k_1-p_3 & p_3 \\
 \end{array}
 \right),~(p_3\neq k_1),$  
 \end{center}\smallskip

 \begin{center}$(3)~\left(
 \begin{array}{ccc}
 k_1 & 0 & 0 \\
 \ell_1 & \ell_2 & k_1-\ell_2 \\
 p_1 & k_1-p_3 & p_3 \\
 \end{array}
 \right),~(\ell_2\neq k_1).$
 \end{center} 
 \end{thm}

 \begin{proof} $N$ is Nijenhuis operator on $(\ll, [,])$ if and only if
 \begin{center}
 $\left\{
 \begin{array}{l}
 k_2 = k_3 = 0\\
 (\ell_3-\ell_2+k_1)(\ell_3+\ell_2- k_1)=0\\
 p_2 (\ell_2- k_1)+ \ell_3(k_1-p_3)=0 \\
 (p_2-p_3+k_1)(p_2 +p_3-k_1)=0\\
 \end{array}
 \right.$.
 \end{center}
 \begin{enumerate}[(I)]
 \item $\ell_3=\ell_2-k_1$, then $\left\{
 \begin{array}{l}
 k_2 = k_3 = 0\\
 \ell_3=\ell_2-k_1\\
 (\ell_2- k_1)(p_2+k_1-p_3)=0 \\
 (p_2-p_3+k_1)(p_2 +p_3-k_1)=0\\
 \end{array}
 \right.$. \smallskip
 \begin{enumerate}
 \item [(IA)] $p_2=-k_1+p_3$, then $\left\{
 \begin{array}{l}
 k_2 = k_3 = 0\\
 \ell_3=\ell_2-k_1\\
 p_2=-k_1+p_3\\
 \end{array}
 \right.$. That is to say,\smallskip
 \begin{center}$\left\{
 \begin{array}{l}
 N(e)=k_1 e\\
 N(f)=\ell_1 e+\ell_2 f+(\ell_2-k_1)g\\
 N(g)=p_1 e+(p_3-k_1) f+p_3 g\\
 \end{array}
 \right.$.
 \end{center}\smallskip
 \item [(IB)] $p_2\neq -k_1+p_3$, we have
 $\left\{
 \begin{array}{l}
 k_2 = k_3 = \ell_3 = 0\\
 \ell_2= k_1 \\
 p_2=-p_3+k_1\\
 k_1\neq p_3\\
 \end{array}
 \right.$. That is to say,\smallskip
 \begin{center}$\left\{
 \begin{array}{l}
 N(e)=k_1 e\\
 N(f)=\ell_1 e+k_1 f\\
 N(g)=p_1 e+(k_1-p_3) f+p_3 g\\
 \end{array}
 \right.~(k_1\neq p_3)$.
 \end{center}
 \end{enumerate}\smallskip
 \item $\ell_3\neq \ell_2-k_1$, then $\left\{
 \begin{array}{l}
 k_2 = k_3 = 0\\
 \ell_3= -\ell_2+ k_1\\
 (\ell_2- k_1)(p_2-k_1+p_3)=0 \\
 (p_2-p_3+k_1)(p_2 +p_3-k_1)=0\\
 \ell_3\neq \ell_2-k_1\\
 \end{array}
 \right.$.  then
 $\left\{
 \begin{array}{l}
 k_2 = k_3 = 0\\
 \ell_3= -\ell_2+ k_1\\
 p_2=k_1-p_3 \\
 \ell_2\neq k_1\\
 \end{array}
 \right.$. That is to say,\smallskip
 \begin{center}$\left\{
 \begin{array}{l}
 N(e)=k_1 e\\
 N(f)=\ell_1 e+\ell_2 f+( k_1-\ell_2) g\\
 N(g)=p_1 e+(k_1-p_3) f+p_3 g\\
 \end{array}
 \right.~(\ell_2\neq k_1)$.
 \end{center}
 \end{enumerate}
 \end{proof}

 \begin{thm}\label{thm:3nb7}  Let $(\ll, [,])$ be a Leibniz algebra and $[,]$ is given by
 \begin{center}
 \begin{tabular}{r|rrr}
 $[,]$ & $e$  & $f$ & $g$ \\
 \hline
 $e$ & $0$  & $0$ & $0$ \\
 $f$ & $0$  & $e$ & $0$  \\
 $g$ & $0$  & $e$ & $\a e$  \\
 \end{tabular},
 where $0\neq \a\in \mathfrak{R}$.
 \end{center}\vskip1mm
 Then all the Nijenhuis operators on $(\ll, [,])$ are given as follows:
 \begin{center}$(1)~\left(
 \begin{array}{ccc}
 k_1 & 0 & 0 \\
 \ell_1 & k_1 & 0 \\
 p_1 & 0 & k_1 \\
 \end{array}
 \right),~(\a > \frac{1}{4}),$ \quad
 $(2)~\left(
 \begin{array}{ccc}
 k_1 & 0 & 0 \\
 \ell_1 & k_1 & 0 \\
 p_1 & \frac{-1-\sqrt{1-4\a}}{2}(p_3- k_1) & p_3 \\
 \end{array}
 \right),~(\a \le \frac{1}{4}, \a \neq 0),$\quad
 \end{center}

 \begin{center}
 $(3)\left(
 \begin{array}{ccc}
 k_1 & 0 & 0 \\
 \ell_1 & k_1 & 0 \\
 p_1 & \frac{-1+\sqrt{1-4\a}}{2}(p_3- k_1) & p_3 \\
 \end{array}
 \right) ~(\a < \frac{1}{4}, \a \neq 0, p_3\neq k_1$), \quad
 \end{center}

 \begin{center}
 $(4)~\left(
 \begin{array}{ccc}
 \ell_2+\frac{1+\sqrt{1-4\a}}{2}\ell_3 & 0 & 0 \\
 \ell_1 & \ell_2 & \ell_3 \\
 p_1 & \frac{1+\sqrt{1-4\a}}{2}(\ell_2+\frac{1+\sqrt{1-4\a}}{2}\ell_3-p_3) & p_3 \\
 \end{array}
 \right),~(\a \le \frac{1}{4}, \a \neq 0, \ell_3\neq 0),$  \quad
 \end{center}

 \begin{center}
 $(5)~\left(
 \begin{array}{ccc}
 \ell_2+\frac{1-\sqrt{1-4\a}}{2}\ell_3 & 0 & 0 \\
 \ell_1 & \ell_2 & \ell_3 \\
 p_1 & \frac{1-\sqrt{1-4\a}}{2}(\ell_2+\frac{1-\sqrt{1-4\a}}{2}\ell_3-p_3) & p_3 \\
 \end{array}
 \right),~(\a < \frac{1}{4}, \a \neq 0,\ell_3\neq 0).$
 \end{center} 
 \end{thm}

 \begin{proof} $N$ is Nijenhuis operator on $(\ll, [,])$ if and only if
 \begin{center}
 $\left\{
 \begin{array}{l}
 k_2 = k_3 = 0\\
 (\ell_2-k_1)^2+\ell_3(\ell_2- k_1)+\a \ell_3^2=0\\
 \a \ell_3 (p_3- k_1)+ p_2 (\ell_2+\ell_3- k_1)=0 \\
 (\ell_2- k_1)(p_3-k_1)-\ell_3 p_2=0\\
 \a (p_3-k_1)^2+p_2(p_3- k_1)+ p_2^2=0\\
 \a \neq 0\\
 \end{array}
 \right.$.
 \end{center}
 \begin{enumerate}[(I)]
 \item $\ell_3=0$, then $\left\{
 \begin{array}{l}
 k_2 = k_3 = \ell_3 = 0\\
 \ell_2=k_1\\
 (\frac{4\a-1}{4}) (p_3-k_1)^2+(p_2+\frac{1}{2}(p_3- k_1))^2=0\\
 \a \neq 0\\
 \end{array}
 \right.$. \smallskip
 \begin{enumerate}
 \item [(IA)] $\a > \frac{1}{4}$, then $\left\{
 \begin{array}{l}
 k_2 = k_3 = \ell_3 = p_2 = 0\\
 p_3=\ell_2=k_1\\
 \end{array}
 \right.$. That is to say,\smallskip
 \begin{center}$\left\{
 \begin{array}{l}
 N(e)=k_1 e\\
 N(f)=\ell_1 e+k_1 f\\
 N(g)=p_1 e+k_1 g\\
 \end{array}
 \right.~(\a > \frac{1}{4})$.
 \end{center}\smallskip

 %\item [(IB)] $\a = \frac{1}{4}$, we have
 %$\left\{
 %\begin{array}{l}
 %k_2 = k_3 = \ell_3 = 0\\
 %\ell_2=k_1\\
 %p_2=\frac{1}{2}(k_1-p_3)\\
 %\end{array}
 %\right.$. That is to say,\smallskip
 %\begin{center}$\left\{
 %\begin{array}{l}
 %N(e)=k_1 e\\
 %N(f)=\ell_1 e+k_1 f\\
 %N(g)=p_1 e+\frac{1}{2}(k_1-p_3) f+p_3 g\\
 %\end{array}
 %\right.~(\a = \frac{1}{4})$.
 %\end{center}

 \item [(IB)] $\a \le \frac{1}{4}$, we have
 $\left\{
 \begin{array}{l}
 k_2 = k_3 = \ell_3 = 0\\
 \ell_2=k_1\\
 (p_2+\frac{1+\sqrt{1-4\a}}{2}(p_3- k_1))(p_2+\frac{1-\sqrt{1-4\a}}{2}(p_3- k_1))=0\\
 \a \neq 0, \a \le \frac{1}{4}\\
 \end{array}
 \right.$.

 \begin{enumerate}
 \item[(IB1)] $p_2=\frac{-1-\sqrt{1-4\a}}{2}(p_3- k_1)$, then $\left\{
 \begin{array}{l}
 k_2 = k_3 = \ell_3 = 0\\
 \ell_2=k_1\\
 p_2=\frac{-1-\sqrt{1-4\a}}{2}(p_3- k_1)\\
 \a \neq 0\\
 \end{array}
 \right.$. That is to say,\smallskip
 \begin{center}$\left\{
 \begin{array}{l}
 N(e)=k_1 e\\
 N(f)=\ell_1 e+k_1 f\\
 N(g)=p_1 e+\frac{-1-\sqrt{1-4\a}}{2}(p_3- k_1) f+p_3 g\\
 \end{array}
 \right.~(\a \le \frac{1}{4}, \a \neq 0)$.
 \end{center}\smallskip

 \item[(IB2)] $p_2\neq \frac{-1-\sqrt{1-4\a}}{2}(p_3- k_1)$, then
 $\left\{
 \begin{array}{l}
 k_2 = k_3 = \ell_3 = 0\\
 \ell_2=k_1\\
 p_2=\frac{-1+\sqrt{1-4\a}}{2}(p_3- k_1)\\
 \a \neq 0, \a < \frac{1}{4}, p_3\neq k_1\\
 \end{array}
 \right.$. That is to say,\smallskip
 \begin{center}$\left\{
 \begin{array}{l}
 N(e)=k_1 e\\
 N(f)=\ell_1 e+k_1 f\\
 N(g)=p_1 e+\frac{-1+\sqrt{1-4\a}}{2}(p_3- k_1) f+p_3 g\\
 \end{array}
 \right.~(\a < \frac{1}{4}, \a \neq 0, p_3\neq k_1)$.
 \end{center}
 \end{enumerate}
 \end{enumerate}

 \item $\ell_3\neq 0$, then $\left\{
 \begin{array}{l}
 k_2 = k_3 = 0\\
 (\ell_2-k_1+\frac{1}{2}\ell_3)^2+\frac{4\a-1}{4} \ell_3^2=0\\
 p_2=\frac{(\ell_2- k_1)(p_3-k_1)}{\ell_3} \\
 \a \neq 0, \ell_3\neq 0\\
 \end{array}
 \right.$. \smallskip
 \begin{enumerate}
 \item [(IIA)] $\a > \frac{1}{4}$, this assumption is not valid since
 $\ell_3\neq 0$.
 %\item [(IIB)] $\a = \frac{1}{4}$, then $\left\{
 %\begin{array}{l}
 %k_2 = k_3 = 0\\
 %\ell_3=2k_1-2\ell_2\\
 %p_2=\frac{1}{2}(k_1 -p_3)
 %k_1\neq \ell_2
 %\end{array}
 %\right.$. That is to say,\smallskip
 %\begin{center}$\left\{
 %\begin{array}{l}
 %N(e)=k_1 e\\
 %N(f)=\ell_1 e+\ell_2 f+ (2k_1-2\ell_2) g\\
 %N(g)=p_1 e+\frac{1}{2}(k_1 -p_3) f+p_3 g\\
 %\end{array}
 %\right.~(\a = \frac{1}{4},  k_1\neq \ell_2)$.
 %\end{center}

 \item[(IIB)] $\a \le \frac{1}{4}$, then $\left\{
 \begin{array}{l}
 k_2 = k_3 = 0\\
 (\ell_2-k_1+\frac{1+\sqrt{1-4\a}}{2}\ell_3)(\ell_2-k_1+
 \frac{1-\sqrt{1-4\a}}{2}\ell_3)=0\\
 p_2=\frac{(\ell_2- k_1)(p_3-k_1)}{\ell_3} \\
 \a \neq 0, \ell_3\neq 0, \a \le \frac{1}{4}\\
 \end{array}
 \right.$.\smallskip
 \begin{enumerate}
 \item [(IIB1)]
 $k_1=\ell_2+\frac{1+\sqrt{1-4\a}}{2}\ell_3$,then
 $\left\{
 \begin{array}{l}
 k_2 = k_3 = 0\\
 k_1=\ell_2+\frac{1+\sqrt{1-4\a}}{2}\ell_3\\
 p_2=\frac{1+\sqrt{1-4\a}}{2}(\ell_2+\frac{1+\sqrt{1-4\a}}{2}\ell_3-p_3)\\
  \a \neq 0, \ell_3\neq 0, \a \le \frac{1}{4}\\
 \end{array}
 \right.$. That is to say,\smallskip
 \begin{center}$\left\{
 \begin{array}{l}
 N(e)=(\ell_2+\frac{1+\sqrt{1-4\a}}{2}\ell_3) e\\
 N(f)=\ell_1 e+\ell_2 f+ \ell_3 g\\
 N(g)=p_1 e+\frac{1+\sqrt{1-4\a}}{2}(\ell_2+\frac{1+\sqrt{1-4\a}}{2}\ell_3-p_3) f+p_3 g\\
 \end{array}
 \right.$.
 $~(\a \le \frac{1}{4}, \ell_3\neq 0, \a \neq 0).$
 \end{center}\smallskip
  \item [(IIC2)] $k_1\neq \ell_2+\frac{1+\sqrt{1-4\a}}{2}\ell_3$, then     $\left\{
 \begin{array}{l}
 k_2 = k_3 = 0\\
 k_1=\ell_2+\frac{1-\sqrt{1-4\a}}{2}\ell_3\\
 p_2=\frac{1-\sqrt{1-4\a}}{2}(\ell_2+\frac{1-\sqrt{1-4\a}}{2}\ell_3-p_3)\\
 \a \neq 0, \a < \frac{1}{4}, \ell_3\neq 0\\
 \end{array}
 \right.$. That is to say,\smallskip
 \begin{center}$\left\{
 \begin{array}{l}
 N(e)=(\ell_2+\frac{1-\sqrt{1-4\a}}{2}\ell_3) e\\
 N(f)=\ell_1 e+\ell_2 f+ \ell_3 g\\
 N(g)=p_1 e+\frac{1-\sqrt{1-4\a}}{2}(\ell_2+\frac{1-\sqrt{1-4\a}}{2}\ell_3-p_3) f+p_3 g\\
 \end{array}
 \right.$.
 $(\a < \frac{1}{4}, \ell_3\neq 0, \a \neq 0,)$
 \end{center}
  \end{enumerate}
  \end{enumerate}
  \end{enumerate}
 \end{proof}

 \begin{thm}\label{thm:3nb8}  Let $(\ll, [,])$ be a Leibniz algebra and $[,]$ is given by
 \begin{center}
 \begin{tabular}{r|rrr}
 $[,]$ & $e$  & $f$ & $g$ \\
 \hline
 $e$ & $0$  & $0$ & $0$ \\
 $f$ & $0$  & $0$ & $0$  \\
 $g$ & $0$  & $e$ & $0$  \\
 \end{tabular}.
 \end{center}\vskip1mm
 Then all the Nijenhuis operators on $(\ll, [,])$ are given as follows:
 \begin{center}$(1)~\left(
 \begin{array}{ccc}
 k_1 & 0 & 0 \\
 \ell_1 & \ell_2 & 0 \\
 p_1 & p_2& k_1 \\
 \end{array}
 \right),~(p_2\neq 0),$ \quad
 $(2)~\left(
 \begin{array}{ccc}
 k_1 & 0 & 0 \\
 \ell_1 & \ell_2 & 0\\
 p_1 & 0 & k_1 \\
 \end{array}
 \right),~(\ell_2\neq k_1),$  \quad
 \end{center}

 \begin{center}
 $(3)~\left(
 \begin{array}{ccc}
 k_1 & 0 & 0 \\
 \ell_1 & k_1 & \ell_3 \\
 p_1 & 0  & p_3\\
 \end{array}
 \right).$
 \end{center} 
 \end{thm}

 \begin{proof} $N$ is Nijenhuis operator on $(\ll, [,])$ if and only if
 \begin{center}
 $\left\{
 \begin{array}{l}
 k_2 = k_3 = 0\\
 \ell_3(\ell_2-k_1)=0\\
 (p_3-k_1)( \ell_2 - k_1)=0 \\
 p_2 \ell_3=0\\
 p_2(p_3- k_1)=0\\
 \end{array}
 \right.$.
 \end{center}
 \begin{enumerate}[(I)]
 \item $p_2=0$, then $\left\{
 \begin{array}{l}
 k_2 = k_3 = p_2 = 0\\
 \ell_3(\ell_2-k_1)=0\\
 (p_3-k_1)( \ell_2 - k_1)=0 \\
 \end{array}
 \right.$. \smallskip
 \begin{enumerate}
 \item [(IA)] $\ell_2=k_1$, then $\left\{
 \begin{array}{l}
 k_2 = k_3 = p_2 = 0\\
 \ell_2 = k_1\\
 \end{array}
 \right.$. That is to say,\smallskip
 \begin{center}$\left\{
 \begin{array}{l}
 N(e)=k_1 e\\
 N(f)=\ell_1 e+k_1 f+\ell_3 g\\
 N(g)=p_1 e+p_3 g\\
 \end{array}
 \right.$.
 \end{center}\smallskip
 \item [(IB)] $\ell_2\neq k_1$, we have
 $\left\{
 \begin{array}{l}
 k_2 = k_3 = p_2 = \ell_3 = 0\\
 p_3 = k_1\\
 \ell_2\neq k_1\\
 \end{array}
 \right.$. That is to say,\smallskip
 \begin{center}$\left\{
 \begin{array}{l}
 N(e)=k_1 e\\
 N(f)=\ell_1 e+\ell_2 f\\
 N(g)=p_1 e+k_1 g\\
 \end{array}
 \right.~(\ell_2\neq k_1)$.
 \end{center}
 \end{enumerate}\smallskip
 \item $p_2\neq 0$, then $\left\{
 \begin{array}{l}
 k_2 = k_3 = \ell_3 = 0\\
 p_2\neq 0\\
 p_3 = k_1\\
 \end{array}
 \right.$. That is to say,\smallskip
 \begin{center}$\left\{
 \begin{array}{l}
 N(e)=k_1 e\\
 N(f)=\ell_1 e+\ell_2 f\\
 N(g)=p_1 e+p_2 f+k_1 g\\
 \end{array}
 \right.~(p_2 \neq 0)$.
 \end{center}
 \end{enumerate}
 \end{proof}

 \begin{thm}\label{thm:3nb9}  Let $(\ll, [,])$ be a Leibniz algebra and $[,]$ is given by
 \begin{center}
 \begin{tabular}{r|rrr}
 $[,]$ & $e$  & $f$ & $g$ \\
 \hline
 $e$ & $0$  & $0$ & $0$ \\
 $f$ & $0$  & $0$ & $0$  \\
 $g$ & $f$  & $e$ & $0$  \\
 \end{tabular}.
 \end{center}\vskip1mm
 Then all the Nijenhuis operators on $(\ll, [,])$ are given as follows:
 \begin{center}$(1)~\left(
 \begin{array}{ccc}
 k_1 & 0 & 0 \\
 0 & k_1 & 0 \\
 p_1 & p_2 & k_1 \\
 \end{array}
 \right),$ \quad
 $(2)~\left(
 \begin{array}{ccc}
 k_1 & 0 & 0 \\
 \ell_1 & k_1-\ell_1 & 0 \\
 p_1 & -p_1 & k_1-\ell_1  \\
 \end{array}
 \right),~(\ell_1 \neq 0),$  \quad
 \end{center}

 \begin{center}$(3)~\left(
 \begin{array}{ccc}
 k_1 & 0 & 0 \\
 \ell_1 & k_1+\ell_1 & 0 \\
 p_1 & p_1 & k_1+\ell_1 \\
 \end{array}
 \right),~(\ell_1 \neq 0)$ \quad
 $(4)~\left(
 \begin{array}{ccc}
 k_1 & 0 & 0 \\
 0 & k_1 & 0 \\
 0 & 0 & p_3 \\
 \end{array}
 \right),~(p_3\neq k_1),$  \quad
 \end{center}

 \begin{center}
 $(5)~\left(
 \begin{array}{ccc}
 k_1 & k_2 & 0 \\
 k_2 & k_1 & 0 \\
 0 & 0  &  p_3 \\
 \end{array}
 \right),~(k_2\neq 0),$\quad
 $(6)\left(
 \begin{array}{ccc}
 k_1 & k_2 & 0 \\
 k_2 & k_1 & 0 \\
 p_1  & p_1  &  k_1 + k_2 \\
 \end{array}
 \right), ~(k_2\neq 0,p_1\neq 0),$ \quad
    \end{center}

 \begin{center}
 $(7)\left(
 \begin{array}{ccc}
 k_1 & k_2 & 0 \\
 k_2 & k_1 & 0 \\
 p_1  & -p_1  &  k_1 - k_2 \\
 \end{array}
 \right), ~(k_2\neq 0,p_1\neq 0),$ \quad
  $(8)~\left(
 \begin{array}{ccc}
 k_1 & k_2 & 0 \\
 k_1 - \ell_2+k_2 & \ell_2 & 0 \\
 p_1 & -p_1  &  \ell_2-k_2 \\
 \end{array}
 \right),~(k_2\neq 0, k_1\neq \ell_2),$\quad
 \end{center}

 \begin{center}
 $(9)~\left(
 \begin{array}{ccc}
 k_1 & k_2 & 0 \\
 \ell_2-k_1+k_2 & \ell_2 & 0 \\
 p_1 & p_1  &  k_1+\ell_2 \\
 \end{array}
 \right),~(k_2\neq 0, k_1\neq \ell_2).$
 \end{center} 
 \end{thm}

 \begin{proof} $N$ is Nijenhuis operator on $(\ll, [,])$ if and only if
 \begin{center}
 $\left\{
 \begin{array}{l}
 \ell_3 = k_3 =0\\
 (k_1 - \ell_2)(p_3 - \ell_2)-k_2(k_2 -\ell_1)=0\\
 k_2(p_3-k_1)-\ell_1(p_3- \ell_2)=0\\
 \ell_1( k_2 - \ell_1)- (k_1-\ell_2)(p_3-k_1)=0 \\
 p_1(p_3-\ell_2)- p_2 k_2=0\\
 p_2(p_3- k_1)- p_1 \ell_1=0\\
 \end{array}
 \right.$.
 \end{center}
 \begin{enumerate}[(I)]
 \item $k_2=0$, then $\left\{
 \begin{array}{l}
 \ell_3 = k_3 = k_2 = 0\\
 (k_1 - \ell_2)(p_3 - \ell_2)=0\\
 \ell_1(p_3- \ell_2)=0\\
 -\ell_1^2- (k_1-\ell_2)(p_3-k_1)=0 \\
 p_1(p_3-\ell_2)=0\\
 p_2(p_3- k_1)- p_1 \ell_1=0\\
 \end{array}
 \right.$. \smallskip
 \begin{enumerate}
 \item [(IA)] $p_3=\ell_2$, then $\left\{
 \begin{array}{l}
 \ell_3 = k_3 = k_2 = 0\\
 (\ell_2-k_1+\ell_1)(\ell_2-k_1-\ell_1)=0 \\
 p_2(\ell_2- k_1)- p_1 \ell_1=0\\
 p_3=\ell_2\\
 \end{array}
 \right.$. \smallskip
 \begin{enumerate}
 \item [(IA1)] $\ell_2=k_1-\ell_1$, then $\left\{
 \begin{array}{l}
 \ell_3 = k_3 = k_2 = 0\\
 \ell_2=k_1-\ell_1 \\
 \ell_1(p_2+ p_1)=0\\
 p_3=\ell_2\\
 \end{array}
 \right.$. \smallskip
 \begin{enumerate}
 \item [(IA1a)]  $\ell_1=0$, then $\left\{
 \begin{array}{l}
 \ell_3 = k_3 = k_2 = \ell_1 = 0\\
 p_3=\ell_2=k_1\\
 \end{array}
 \right.$.  That is to say,\smallskip
 \begin{center}$\left\{
 \begin{array}{l}
 N(e)=k_1 e\\
 N(f)=k_1 f \\
 N(g)=p_1 e+p_2 f+k_1 g\\
 \end{array}
 \right. $.
 \end{center}\smallskip
 \item [(IA1b)] $\ell_1 \neq 0$, we have $\left\{
 \begin{array}{l}
 \ell_3 = k_3 = k_2 = 0\\
 \ell_2=k_1-\ell_1 \\
 p_2=- p_1\\
 p_3=k_1-\ell_1\\
 \ell_1 \neq 0
 \end{array}
 \right.$. That is to say,\smallskip
 \begin{center}$\left\{
 \begin{array}{l}
 N(e)=k_1 e\\
 N(f)=\ell_1 e+(k_1-\ell_1) f \\
 N(g)=p_1 e-p_1 f+(k_1-\ell_1) g\\
 \end{array}
 \right. ~(\ell_1 \neq 0)$.
 \end{center}
 \end{enumerate}\smallskip
 \item [(IA2)] $\ell_2\neq k_1-\ell_1$, we have
 $\left\{
 \begin{array}{l}
 \ell_3 = k_3 = k_2 = 0\\
 \ell_2=k_1+\ell_1\\
 p_2= p_1\\
 p_3=k_1+\ell_1\\
 \ell_1\neq 0\\
 \end{array}
 \right.$. That is to say,\smallskip
 \begin{center}$\left\{
 \begin{array}{l}
 N(e)=k_1 e\\
 N(f)=\ell_1 e+(k_1+\ell_1) f \\
 N(g)=p_1 e+p_1 f+(k_1+\ell_1) g\\
 \end{array}
 \right. ~(\ell_1 \neq 0)$.
 \end{center}
 \end{enumerate}\smallskip

 \item [(IB)]  $p_3\neq \ell_2$, we have
 $\left\{
 \begin{array}{l}
 \ell_3 = k_3 = k_2 = p_1 = \ell_1 = p_2 = 0\\
 \ell_2=k_1\\
 p_3\neq k_1\\
 \end{array}
 \right.$.  That is to say,\smallskip
 \begin{center}$\left\{
 \begin{array}{l}
 N(e)=k_1 e\\
 N(f)=k_1 f \\
 N(g)=p_3 g\\
 \end{array}
 \right. ~(p_3\neq k_1)$.
 \end{center}
 \end{enumerate}

 \item $k_2 \neq0$, then $\left\{
 \begin{array}{l}
 \ell_3 = k_3 =0\\
 \ell_1=\frac{(k_1 - \ell_2)(\ell_2-p_3 )}{k_2}+k_2 \\
 (k_1- \ell_2)(k_2-\ell_2+p_3)(k_2+\ell_2-p_3)=0\\
 p_1(p_3-\ell_2)- p_2 k_2=0\\
 p_2(p_3- k_1)- p_1 \ell_1=0\\
 k_2 \neq 0\\
 \end{array}
 \right.$. \smallskip
 \begin{enumerate}
 \item [(IIA)] $k_1=\ell_2$, then
 $\left\{
 \begin{array}{l}
 \ell_3 = k_3 =0\\
 \ell_1=k_2 \\
 k_1=\ell_2\\
 p_2(p_3-\ell_2-\ell_1)(p_3-\ell_2+\ell_1) =0\\
 p_1=\frac{p_2(p_3- \ell_2)}{\ell_1}\\
 k_2 \neq 0\\
 \end{array}
 \right.$.\smallskip
 \begin{enumerate}
 \item [(IIA1)] $p_2=0$, then $\left\{
 \begin{array}{l}
 \ell_3 = k_3 =p_1 = p_2 = 0\\
 \ell_1=k_2 \\
 \ell_2=k_1\\
 k_2 \neq 0\\
 \end{array}
 \right.$. That is to say, \smallskip
 \begin{center}$\left\{
 \begin{array}{l}
 N(e)=k_1 e+k_2 f\\
 N(f)=k_2 e+k_1 f\\
 N(g)=p_3 g\\
 \end{array}
 \right.~(k_2 \neq 0)$.
 \end{center}\smallskip
 \item [(IIA2)] $p_2\neq 0$, then $\left\{
 \begin{array}{l}
 \ell_3 = k_3 =0\\
 \ell_1=k_2 \\
 k_1=\ell_2\\
(p_3-\ell_2-\ell_1)(p_3-\ell_2+\ell_1) =0\\
 p_1=\frac{p_2(p_3- \ell_2)}{\ell_1}\\
 k_2 \neq 0, p_2\neq 0\\
 \end{array}
 \right.$.\smallskip
 \begin{enumerate}
 \item [(IIA2a)] $p_3=\ell_2+\ell_1$, then $\left\{
 \begin{array}{l}
 \ell_3 = k_3 =0\\
 \ell_1=k_2 \\
 k_1=\ell_2\\
 p_3=k_1+k_2\\
 p_1=p_2\\
 k_2 \neq 0, p_1\neq 0\\
 \end{array}
 \right.$.
 That is to say,\smallskip
 \begin{center}$\left\{
 \begin{array}{l}
 N(e)=k_1 e+k_2 f\\
 N(f)=k_2 e+k_1 f \\
 N(g)=p_1 e+p_1 f+(k_1+k_2) g\\
 \end{array}
 \right.~( k_2 \neq 0, p_1\neq 0)$.
 \end{center}\smallskip
  \item [(IIA2b)] $p_3\neq \ell_2+\ell_1$, then $\left\{
 \begin{array}{l}
 \ell_3 = k_3 =0\\
 \ell_1=k_2 \\
 k_1=\ell_2\\
 p_3=k_1-k_2\\
 p_1=-p_2\\
 k_2 \neq 0, p_1\neq 0\\
 \end{array}
 \right.$.
 That is to say,\smallskip
 \begin{center}$\left\{
 \begin{array}{l}
 N(e)=k_1 e+k_2 f\\
 N(f)=k_2 e+k_1 f \\
 N(g)=p_1 e-p_1 f+(k_1-k_2) g\\
 \end{array}
 \right.~( k_2 \neq 0, p_1\neq 0)$.
 \end{center}
 \end{enumerate}
 \end{enumerate}\smallskip

 \item [(IIB)] $k_1\neq \ell_2$, then
 $\left\{
 \begin{array}{l}
 \ell_3 = k_3 =0\\
 \ell_1=\frac{(k_1 - \ell_2)(\ell_2-p_3 )}{k_2}+k_2 \\
 (k_2-\ell_2+p_3)(k_2+\ell_2-p_3)=0\\
 p_1(p_3-\ell_2)- p_2 k_2=0\\
 p_2(p_3- k_1)- p_1 \ell_1=0\\
 k_2 \neq 0, k_1\neq \ell_2\\
 \end{array}
 \right.$.\smallskip
 \begin{enumerate}
 \item [(IIB1)] $p_3=\ell_2-k_2$, then $\left\{
 \begin{array}{l}
 \ell_3 = k_3 =0\\
 \ell_1=k_1 - \ell_2+k_2 \\
 p_1 =-p_2\\
 p_3=\ell_2-k_2\\
 k_2 \neq 0, k_1\neq \ell_2\\
 \end{array}
 \right.$. That is to say, \smallskip
 \begin{center}$\left\{
 \begin{array}{l}
 N(e)=k_1 e+k_2 f\\
 N(f)=(k_1 - \ell_2+k_2) e+\ell_2 f\\
 N(g)=p_1 e-p_1 f +(\ell_2-k_2) g\\
 \end{array}
 \right.~(k_2 \neq 0,k_1\neq \ell_2)$.
 \end{center}\smallskip
 \item [(IIB2)] $p_3\neq \ell_2-k_2$, then $\left\{
 \begin{array}{l}
 \ell_3 = k_3 =0\\
 \ell_1= \ell_2-k_1+k_2 \\
 p_3 =k_2+\ell_2\\
 p_1= p_2\\
 k_2 \neq 0,k_1\neq \ell_2\\
 \end{array}
 \right.$.
 That is to say,\smallskip
 \begin{center}$\left\{
 \begin{array}{l}
 N(e)=k_1 e+k_2 f\\
 N(f)=(\ell_2-k_1+k_2) e+\ell_2 f \\
 N(g)=p_1 e+p_1 f+(k_2+\ell_2) g\\
 \end{array}
 \right.~(  k_2 \neq 0,k_1\neq \ell_2)$.
 \end{center}
 \end{enumerate}
 \end{enumerate}
 \end{enumerate}
 \end{proof}

 \begin{thm}  \label{thm:3nb10}  Let $(\ll, [,])$ be a Leibniz algebra and $[,]$ is given by
 \begin{center}
 \begin{tabular}{r|rrr}
 $[,]$ & $e$  & $f$ & $g$ \\
 \hline
 $e$ & $0$  & $0$ & $0$ \\
 $f$ & $0$  & $0$ & $0$  \\
 $g$ & $f$  & $-e$ & $0$  \\
 \end{tabular}.
 \end{center}\vskip1mm
 Then all the Nijenhuis operators on $(\ll, [,])$ are given as follows:
 \begin{center}$(1)~\left(
 \begin{array}{ccc}
 k_1 & 0 & 0 \\
 0 & k_1 & 0 \\
 p_1 & p_2 & k_1 \\
 \end{array}
 \right),$ \quad
 $(2)\left(
 \begin{array}{ccc}
 k_1 & 0 & 0 \\
 0 & k_1 & 0 \\
 0 & 0 & p_3 \\
 \end{array}
 \right) ~(p_3\neq k_1),$ \quad
 \end{center}

  \begin{center}
 $(3)~\left(
 \begin{array}{ccc}
 k_1 & k_2 & 0 \\
 -k_2 & k_1 & 0 \\
 0 & 0 & p_3 \\
 \end{array}
 \right),~(k_2\neq 0).$
 \end{center} 
 \end{thm}

 \begin{proof} $N$ is Nijenhuis operator on $(\ll, [,])$ if and only if
 \begin{center}
 $\left\{
 \begin{array}{l}
 k_3= \ell_3 =0\\
 (k_1-\ell_2)(p_3-\ell_2)+k_2(\ell_1+ k_2)=0\\
 \ell_1 (\ell_2- p_3)+ k_2(k_1-p_3)=0 \\
 (k_1-\ell_2)(p_3-k_1)-\ell_1( \ell_1+k_2)=0\\
 p_1(p_3- \ell_2)+ p_2 k_2=0\\
 p_2(k_1-p_3)-p_1 \ell_1=0\\
 \end{array}
 \right.$.
 \end{center}
 \begin{enumerate}[(I)]
 \item $k_2=0$, then $\left\{
 \begin{array}{l}
 k_3= \ell_3 = k_2 = 0\\
 (k_1-\ell_2)(p_3-\ell_2)=0\\
 \ell_1 (\ell_2- p_3)=0 \\
 (k_1-\ell_2)(p_3-k_1)-\ell_1^2=0\\
 p_1(p_3- \ell_2)=0\\
 p_2(k_1-p_3)-p_1 \ell_1=0\\
 \end{array}
 \right.$. \smallskip
 \begin{enumerate}
 \item [(IA)] $p_3=\ell_2$, then $\left\{
 \begin{array}{l}
 k_3= \ell_3 = k_2 = \ell_1= 0\\
 p_3=\ell_2=k_1\\
 \end{array}
 \right.$. \smallskip
  That is to say,\smallskip
 \begin{center}$\left\{
 \begin{array}{l}
 N(e)=k_1 e\\
 N(f)=k_1 f\\
 N(g)=p_1 e+ p_2 f +k_1 g\\
 \end{array}
 \right. $.
 \end{center}

 \item [(IB)] $p_3\neq \ell_2$, we have
 $\left\{
 \begin{array}{l}
 k_3= \ell_3 = k_2 = \ell_1 = p_1 = p_2 = 0\\
 \ell_2=k_1\\
 p_3\neq k_1\\
 \end{array}
 \right.$. That is to say,\smallskip
 \begin{center}$\left\{
 \begin{array}{l}
 N(e)=k_1 e\\
 N(f)=k_1 f\\
 N(g)=p_3 g\\
 \end{array}
 \right.~(p_3\neq k_1)$.
 \end{center}
  \end{enumerate}

 \item $k_2\neq 0$, then $\left\{
 \begin{array}{l}
 k_3= \ell_3 =0\\
 \ell_1= \frac{(k_1-\ell_2)(\ell_2-p_3)}{k_2}-k_2\\
 (k_1-\ell_2)((\ell_2- p_3)^2+ k_2^2)=0 \\
 (k_1-\ell_2)(p_3-k_1)-\ell_1( \ell_1+k_2)=0\\
 p_2= \frac{p_1(\ell_2-p_3)}{k_2}\\
 p_1(\ell_2-p_3)(k_1-p_3)-p_1 \ell_1 k_2=0\\
 k_2\neq 0\\
 \end{array}
 \right.$.

 Then $\left\{
 \begin{array}{l}
 k_3= \ell_3 = p_2= p_1= 0\\
 \ell_1= -k_2\\
 \ell_2=k_1\\
 k_2\neq 0\\
 \end{array}
 \right.$. That is to say,\smallskip
 \begin{center}$\left\{
 \begin{array}{l}
 N(e)=k_1 e+k_2 f\\
 N(f)=-k_2 e+k_1 f\\
 N(g)=p_3 g\\
 \end{array}
 \right.~(k_2\neq 0)$.
 \end{center}\smallskip
 \end{enumerate}
 \end{proof}

 \begin{thm}\label{thm:3nb11}  Let $(\ll, [,])$ be a Leibniz algebra and $[,]$ is given by
 \begin{center}
 \begin{tabular}{r|rrr}
 $[,]$ & $e$  & $f\quad$ & $g$ \\
 \hline
 $e$ & $0$  & $0\quad$ & $0$ \\
 $f$ & $0$  & $0\quad$ & $0$  \\
 $g$ & $f$  & $\a e+f$ & $0$  \\
 \end{tabular},
 where $0\neq \a\in \mathfrak{R}$.
 \end{center}\vskip1mm
 Then all the Nijenhuis operators on $(\ll, [,])$ are given as follows:
\begin{center}$(1)~\left(
 \begin{array}{ccc}
 k_1 & 0 & 0 \\
 0 & k_1 & 0 \\
 p_1 & 0 & k_1 \\
 \end{array}
 \right),~(p_1\neq 0),$ \quad
 $(2)~\left(
 \begin{array}{ccc}
 k_1 & k_2 & 0 \\
 0 & k_1+\frac{1+\sqrt{1+4\a}}{2}k_2 & 0 \\
 0 & 0 & k_1 \\
 \end{array}
 \right),~(k_2\neq 0, \a \ge \frac{-1}{4}, \a\neq 0),$  \quad
 \end{center}

 \begin{center}
 $(3)~\left(
 \begin{array}{ccc}
 k_1 & k_2 & 0 \\
 0 & k_1+\frac{1-\sqrt{1+4\a}}{2}k_2 & 0 \\
 0 & 0 & k_1 \\
 \end{array}
 \right),~(k_2\neq 0, \a > \frac{-1}{4}, \a\neq 0),$\quad
 $(4)\left(
 \begin{array}{ccc}
 k_1 & 0 & 0 \\
 0 & k_1 & 0 \\
 0 & 0 & p_3 \\
 \end{array}
 \right) $ \quad
 \end{center}

   \begin{center}
 $(5)~\left(
 \begin{array}{ccc}
 k_1 & k_2 & 0 \\
 \a k_2 & k_1+ k_2 & 0 \\
 0 & 0 & p_3 \\
 \end{array}
 \right),~( \a \neq 0),$\quad
  \end{center}

 \begin{center} 
  $(6)~\left(
 \begin{array}{ccc}
 k_1& k_2 & 0 \\
 \frac{-1+\sqrt{1+4\a}}{2}(p_3-k_1) & (\frac{1-\sqrt{1+4\a}}{2}k_2+p_3)& 0 \\
 0 & 0 & p_3\\
 \end{array}
 \right), ~$\quad\\
 $(k_2\neq \frac{-1+\sqrt{1+4\a}}{2\a}(p_3-k_1), p_3\neq k_1, \a> \frac{-1}{4}, \a \neq 0),$
 \end{center}\smallskip

 \begin{center} 
 $(7)~\left(
 \begin{array}{ccc}
 k_1& k_2 & 0 \\
 \frac{-1-\sqrt{1+4\a}}{2}(p_3-k_1) & (\frac{1+\sqrt{1+4\a}}{2}k_2+p_3)& 0 \\
 0 & 0 & p_3\\
 \end{array}
 \right),~$\\
 $\quad (k_2\neq \frac{-1-\sqrt{1+4\a}}{2\a}(p_3-k_1), p_3\neq k_1, \a\ge \frac{-1}{4}, \a \neq 0),$
 \end{center}\smallskip

   \begin{center}
  $(8)~\left(
 \begin{array}{ccc}
 k_1 & 0 & 0 \\
 0 & k_1 & 0 \\
 p_1 & p_2 & k_1\\
 \end{array}
 \right),~(p_2\neq 0, \a<\frac{-1}{4}),$\smallskip

 $(9)~\left(
 \begin{array}{ccc}
 k_1 & 0 & 0 \\
 0 & k_1 & 0 \\
 p_1 & p_2 & k_1\\
 \end{array}
 \right),~(p_2\neq 0, \a\ge\frac{-1}{4}, p_1\neq \frac{-1\pm \sqrt{1+4\a}}{2}p_2),$ 
 \end{center}\smallskip

 \begin{center}
 \small{
 $(10)~\left(
 \begin{array}{ccc}
 k_1 & \frac{1-\sqrt{1+4\a}}{2\a}(k_1-\ell_2)+\frac{(1-\sqrt{1+4\a})^2}{4\a^2}\ell_1 & 0 \\
 \ell_1 & \ell_2 & 0 \\
 \frac{-1-\sqrt{1+4\a}}{2}p_2 & p_2 & k_1+\frac{1-\sqrt{1+4\a}}{2\a} \ell_1 \\
 \end{array}
 \right),~(p_2\neq 0,  \a \neq 0, \a \ge \frac{-1}{4})$\smallskip

  $(11)~\left(
 \begin{array}{ccc}
 k_1 & \frac{1+\sqrt{1+4\a}}{2\a}(k_1-\ell_2)+\frac{(1+\sqrt{1+4\a})^2}{4\a^2}\ell_1 & 0 \\
 \ell_1 & \ell_2 & 0 \\
 \frac{-1+\sqrt{1+4\a}}{2}p_2 & p_2 & k_1+\frac{1+\sqrt{1+4\a}}{2\a} \ell_1 \\
 \end{array}
 \right),~(p_2\neq 0, \a \neq 0, \a > \frac{-1}{4}).$}
 \end{center} 
 \end{thm}

 \begin{proof} $N$ is Nijenhuis operator on $(\ll, [,])$ if and only if
 \begin{center}
 $\left\{
 \begin{array}{l}
 \ell_3 = k_3 =0\\
 (p_3- \ell_2)(k_1+k_2-\ell_2)+k_2(\ell_1- \a k_2)=0\\
 \a k_2( p_3- k_1)+ \ell_1(\ell_2-p_3-k_2)=0 \\
 \a k_1(k_1-p_3-\ell_2)+\ell_1(k_1-p_3-\ell_1)+\a p_3 \ell_2+\a k_2 \ell_1=0\\
 (p_3-\ell_2)(p_1+ p_2)-\a p_2 k_2=0\\
 \a p_2(p_3-k_1)-\ell_1(p_1+ p_2)=0\\
 \end{array}
 \right.$.
 \end{center}
 \begin{enumerate}[(I)]
 \item $p_2=0$, then $\left\{
 \begin{array}{l}
 \ell_3 = k_3 = p_2 = 0\\
 (p_3- \ell_2)(k_1+k_2-\ell_2)+k_2(\ell_1- \a k_2)=0\\
 \a k_2( p_3- k_1)+ \ell_1(\ell_2-p_3-k_2)=0 \\
 \a k_1(k_1-p_3-\ell_2)+\ell_1(k_1-p_3-\ell_1)+\a p_3 \ell_2+\a k_2 \ell_1=0\\
 (p_3-\ell_2)p_1=0\\
 \ell_1 p_1=0\\
 \end{array}
 \right.$. \smallskip
 \begin{enumerate}
 \item [(IA)] $p_1\neq 0$, then $\left\{
 \begin{array}{l}
 \ell_3 = k_3 = p_2 = \ell_1= k_2=0\\
 p_3=\ell_2=k_1\\
 \end{array}
 \right.$. That is to say,\smallskip
 \begin{center}$\left\{
 \begin{array}{l}
 N(e)=k_1 e\\
 N(f)=k_1 f\\
 N(g)=p_1 e+k_1 g\\
 \end{array}
 \right.~(p_1\neq 0)$.
 \end{center}\smallskip

 \item [(IB)] $p_1=0$, then $\left\{
 \begin{array}{l}
 \ell_3 = k_3 = p_2 = p_1=0\\
 (p_3- \ell_2)(k_1+k_2-\ell_2)+k_2(\ell_1- \a k_2)=0\\
 \a k_2( p_3- k_1)+ \ell_1(\ell_2-p_3-k_2)=0 \\
 \a k_1(k_1-p_3-\ell_2)+\ell_1(k_1-p_3-\ell_1)+\a p_3 \ell_2+\a k_2 \ell_1=0\\
 \end{array}
 \right.$.\smallskip
 \begin{enumerate}
 \item [(IB1)] $\ell_1 = 0$, then $\left\{
 \begin{array}{l}
 \ell_3 = k_3 = p_2 = \ell_1 = p_1=0\\
 (p_3- \ell_2)(k_1+k_2-\ell_2)- \a k_2^2=0\\
 \a k_2( p_3- k_1)=0 \\
 \a k_1(k_1-p_3-\ell_2)+\a p_3 \ell_2=0\\
 \end{array}
 \right.$. \smallskip
 \begin{enumerate}
 \item [(IB1a)] $k_2 = 0$, then $\left\{
 \begin{array}{l}
 \ell_3 = k_3 = p_2 = \ell_1 = p_1= k_2 =0\\
 (p_3- \ell_2)(k_1-\ell_2)=0\\
 \a k_1(k_1-p_3-\ell_2)+\a p_3 \ell_2=0\\
 \end{array}
 \right.$. \smallskip
 \begin{itemize}
 \item [(IB1ai)]  $k_1=\ell_2$, then  $\left\{
 \begin{array}{l}
 \ell_3 = k_3 = p_2 = \ell_1 = p_1= k_2 =0\\
 \ell_2=k_1\\
 \end{array}
 \right.$.  That is to say,\smallskip
 \begin{center}$\left\{
 \begin{array}{l}
 N(e)=k_1 e\\
 N(f)=k_1 f\\
 N(g)=p_3 g\\
 \end{array}
 \right.$.
 \end{center}\smallskip

 \item [(IB1aii)] $k_1\neq \ell_2$, we have $\left\{
 \begin{array}{l}
 \ell_3 = k_3 = p_2 = \ell_1 = p_1= k_2 =0\\
 p_3= \ell_2\\
 \ell_2=k_1\\
 \end{array}
 \right.$. This assumption is not valid since $k_1\neq \ell_2$.
 \end{itemize}\smallskip

 \item [(IB1b)]
 $k_2\neq 0$, then $\left\{
 \begin{array}{l}
 \ell_3 = k_3 = p_2 = \ell_1 = p_1=0\\
 (k_1- \ell_2+\frac{1}{2}k_2)^2- \frac{1+4\a}{4} k_2^2=0\\
  p_3= k_1 \\
 k_2\neq 0\\
 \end{array}
 \right.$.
 \begin{itemize}
 \item [(IB1bi)]
 $ \a<\frac{-1}{4}$, This assumption is not valid since $k_2\neq 0$.
 \item [(IB1bii)]
 $ \a\ge\frac{-1}{4}$,
 then $\left\{
 \begin{array}{l}
 \ell_3 = k_3 = p_2 = \ell_1 = p_1=0\\
 (k_1- \ell_2+\frac{1+\sqrt{1+4\a}}{2}k_2)(k_1- \ell_2+\frac{1-\sqrt{1+4\a}}{2}k_2)=0\\
 p_3= k_1 \\
 k_2\neq 0, \a\ge\frac{-1}{4}, \a\neq 0\\
 \end{array}
 \right.$.\smallskip
 \begin{itemize}
 \item [(IB1bii$\a$)]
 $ \ell_2=k_1+\frac{1+\sqrt{1+4\a}}{2}k_2$,
 then $\left\{
 \begin{array}{l}
 \ell_3 = k_3 = p_2 = \ell_1 = p_1=0\\
 \ell_2=k_1+\frac{1+\sqrt{1+4\a}}{2}k_2\\
  p_3= k_1 \\
 k_2\neq 0, \a\ge\frac{-1}{4}, \a\neq 0\\
 \end{array}
 \right.$.
 That is to say,\smallskip
 \begin{center}$\left\{
 \begin{array}{l}
 N(e)=k_1 e+ k_2 f\\
 N(f)=(k_1+\frac{1+\sqrt{1+4\a}}{2}k_2) f\\
 N(g)= k_1 g\\
 \end{array}
 \right.~(k_2\neq 0, \a\ge\frac{-1}{4}, \a\neq 0)$.
 \end{center}\smallskip
 \item [(IB1bii$\b$)]
 $ \ell_2\neq k_1+\frac{1+\sqrt{1+4\a}}{2}k_2$,
 then $\left\{
 \begin{array}{l}
 \ell_3 = k_3 = p_2 = \ell_1 = p_1=0\\
 \ell_2=k_1+\frac{1-\sqrt{1+4\a}}{2}k_2\\
  p_3= k_1 \\
 k_2\neq 0, \a>\frac{-1}{4}, \a\neq 0\\
 \end{array}
 \right.$.
 That is to say,\smallskip
 \begin{center}$\left\{
 \begin{array}{l}
 N(e)=k_1 e+ k_2 f\\
 N(f)=(k_1+\frac{1-\sqrt{1+4\a}}{2}k_2) f\\
 N(g)= k_1 g\\
 \end{array}
 \right.~(k_2\neq 0, \a>\frac{-1}{4}, \a\neq 0)$.
 \end{center}\smallskip
 \end{itemize}
 \end{itemize}
 \end{enumerate}

 \item [(IB2)]
 $\ell_1\neq 0$, then $\left\{
 \begin{array}{l}
 \ell_3 = k_3 = p_2 = p_1=0\\
 \ell_2=\frac{\a k_2(k_1- p_3)}{\ell_1}+p_3+k_2\\
 k_2 (\a k_2- \ell_1)(\ell^2+\ell(p_3-k_1)-\a (p_3-k_1)^2)=0\\
 (\a k_2- \ell_1)(\ell^2+\ell(p_3-k_1)-\a (p_3-k_1)^2)=0\\
 \ell_1\neq 0\\
 \end{array}
 \right.$.

 then
 $\left\{
 \begin{array}{l}
 \ell_3 = k_3 = p_2 = p_1=0\\
 \ell_2=\frac{\a k_2(k_1- p_3)}{\ell_1}+p_3+k_2\\
 (\a k_2- \ell_1)(\ell_1^2+\ell_1(p_3-k_1)-\a (p_3-k_1)^2)=0\\
 \ell_1\neq 0\\
 \end{array}
 \right.$.
 \begin{enumerate}
 \item [(IB2a)] $\a k_2= \ell_1$, then $\left\{
 \begin{array}{l}
 \ell_3 = k_3 = p_2 = p_1=0\\
 \ell_2= k_1+k_2\\
 \ell_1=\a k_2\\
 \a \neq 0\\
 \end{array}
 \right.$. That is to say,\smallskip
 \begin{center}$\left\{
 \begin{array}{l}
 N(e)=k_1 e+ k_2 f\\
 N(f)=\a k_2 e+ (k_1+k_2) f\\
 N(g)=p_3 g \\
 \end{array}
 \right.~(\a \neq 0)$.
 \end{center}\smallskip
 \item [(IB2b)] $\a k_2\neq \ell_1$,  then
 $\left\{
 \begin{array}{l}
 \ell_3 = k_3 = p_2 = p_1=0\\
 \ell_2=\frac{\a k_2(k_1- p_3)}{\ell_1}+p_3+k_2\\
 (\ell_1+\frac{1+\sqrt{1+4\a}}{2}(p_3-k_1))(\ell_1+\frac{1-\sqrt{1+4\a}}{2}(p_3-k_1))=0\\
 \ell_1\neq 0, \a k_2\neq \ell_1, k_1\neq p_3,\a\ge \frac{-1}{4}, \a \neq 0\\
 \end{array}
 \right.$.\smallskip
 \begin{itemize}
 \item [(IB2bi)]
 $\ell_1=\frac{-1-\sqrt{1+4\a}}{2}(p_3-k_1)$, then $\left\{
 \begin{array}{l}
 \ell_3 = k_3 = p_2 = p_1=0\\
 \ell_2=\frac{1+\sqrt{1+4\a}}{2}k_2+p_3\\
 \ell_1=\frac{-1-\sqrt{1+4\a}}{2}(p_3-k_1)\\
 p_3\neq k_1, k_2\neq \frac{-1-\sqrt{1+4\a}}{2\a}(p_3-k_1)\\
 \a\ge \frac{-1}{4}, \a \neq 0\\
 \end{array}
 \right.$.
 That is to say,\smallskip
 \begin{center}$\left\{
 \begin{array}{l}
 N(e)=k_1 e+ k_2 f\\
 N(f)=\frac{-1-\sqrt{1+4\a}}{2}(p_3-k_1) e+ (\frac{1+\sqrt{1+4\a}}{2}k_2+p_3) f\\
 N(g)=p_3 g \\
 \end{array}
 \right.$.
 $( p_3\neq k_1,  k_2\neq \frac{-1-\sqrt{1+4\a}}{2\a}(p_3-k_1), \a\ge \frac{-1}{4}, \a \neq 0).$
 \end{center}\smallskip
 \item [(IB2bii)]
 $\ell_1\neq \frac{-1-\sqrt{1+4\a}}{2}(p_3-k_1)$, then $\left\{
 \begin{array}{l}
 \ell_3 = k_3 = p_2 = p_1=0\\
 \ell_2=\frac{1-\sqrt{1+4\a}}{2}k_2+p_3\\
 \ell_1=\frac{-1+\sqrt{1+4\a}}{2}(p_3-k_1)\\
  p_3\neq k_1, k_2\neq \frac{-1+\sqrt{1+4\a}}{2\a}(p_3-k_1)\\
  \a> \frac{-1}{4}, \a \neq 0\\
 \end{array}
 \right.$.
 That is to say,\smallskip
 \begin{center}$\left\{
 \begin{array}{l}
 N(e)=k_1 e+ k_2 f\\
 N(f)=\frac{-1+\sqrt{1+4\a}}{2}(p_3-k_1) e+ (\frac{1-\sqrt{1+4\a}}{2}k_2+p_3) f\\
 N(g)=p_3 g
 \end{array}
 \right.$.
 $( p_3\neq k_1, k_2\neq \frac{-1+\sqrt{1+4\a}}{2\a}(p_3-k_1), \a> \frac{-1}{4}, \a \neq 0).$
 \end{center}\smallskip
 \end{itemize}
 \end{enumerate}
 \end{enumerate}
 \end{enumerate}

 \item $p_2\neq 0$, then $\left\{
 \begin{array}{l}
 \ell_3 = k_3 =0\\
 p_2\neq 0\\
 \big(\a p_2(k_1-\ell_2)+\ell_1(p_1+p_2)\big)^2\big((p_1+p_2)^2-p_2(p_1+p_2)-\a p_2^2\big)=0\\
 \ell_1\big(\a p_2(k_1-\ell_2)+\ell_1(p_1+p_2)\big)((p_1+p_2)^2-p_2(p_1+p_2)-\a p_2^2)=0\\
 \ell_1^2((p_1+p_2)^2-p_2(p_1+p_2)-\a p_2^2)=0\\
 k_2=\frac{(k_1-\ell_2)(p_1+ p_2)}{\a p_2}+\frac{\ell_1(p_1+ p_2)^2}{\a^2 p_2^2}\\
 p_3=\frac{\ell_1(p_1+ p_2)}{\a p_2}+ k_1\\
 \end{array}
 \right.$.

 \begin{enumerate}
 \item[(IIA)] $(p_1+p_2)^2-p_2(p_1+p_2)-\a p_2^2=0$, then $\left\{
 \begin{array}{l}
 \ell_3 = k_3 =0\\
 p_2\neq 0, \a \ge \frac{-1}{4}, \a \neq 0\\
 (p_1+ \frac{1+\sqrt{1+4\a}}{2}p_2)(p_1+ \frac{1-\sqrt{1+4\a}}{2}p_2)=0\\
 k_2=\frac{(k_1-\ell_2)(p_1+ p_2)}{\a p_2}+\frac{\ell_1(p_1+ p_2)^2}{\a^2 p_2^2}\\
 p_3=\frac{\ell_1(p_1+ p_2)}{\a p_2}+ k_1\\
 \end{array}
 \right.$.\smallskip

 \begin{enumerate}
 \item[(IIA1)] $ p_1= \frac{-1-\sqrt{1+4\a}}{2}p_2$, then
 $\left\{
 \begin{array}{l}
 \ell_3 = k_3 =0\\
 p_2\neq 0,\a \ge \frac{-1}{4}, \a \neq 0\\
 p_1= \frac{-1-\sqrt{1+4\a}}{2}p_2\\
 k_2=\frac{1-\sqrt{1+4\a}}{2\a}(k_1-\ell_2)+  \frac{(1-\sqrt{1+4\a})^2}{4 \a^2}\ell_1\\
 p_3=\frac{1-\sqrt{1+4\a}}{2\a}\ell_1+ k_1\\
 \end{array}
 \right.$.  That is to say,\smallskip
 \begin{center}\hspace{-4mm}$\left\{
 \begin{array}{l}
 N(e)=k_1 e+ (\frac{1-\sqrt{1+4\a}}{2\a}(k_1-\ell_2)+  \frac{(1-\sqrt{1+4\a})^2}{4 \a^2}\ell_1) f\\
 N(f)=\ell_1 e+ \ell_ 2 f\\
 N(g)=\frac{-1-\sqrt{1+4\a}}{2}p_2 e +p_2 f + (\frac{1-\sqrt{1+4\a}}{2\a}\ell_1+ k_1) g \\
 \end{array}
 \right.~(p_2\neq 0,\a \ge \frac{-1}{4}, \a \neq 0).$ 
 \end{center}\smallskip

 \item[(IIA2)] $ p_1\neq \frac{-1-\sqrt{1+4\a}}{2}p_2$, then
 $\left\{
 \begin{array}{l}
 \ell_3 = k_3 =0\\
 p_2\neq 0, \a > \frac{-1}{4}, \a \neq 0\\
 p_1= \frac{-1+\sqrt{1+4\a}}{2}p_2\\
 k_2=\frac{1+\sqrt{1+4\a}}{2\a}(k_1-\ell_2)+  \frac{(1+\sqrt{1+4\a})^2}{4 \a^2}\ell_1\\
 p_3=\frac{1+\sqrt{1+4\a}}{2\a}\ell_1+ k_1\\
 \end{array}
 \right.$.  That is to say,\smallskip
 \begin{center}$\left\{
 \begin{array}{l}
 N(e)=k_1 e+ (\frac{1+\sqrt{1+4\a}}{2\a}(k_1-\ell_2)+  \frac{(1+\sqrt{1+4\a})^2}{4 \a^2}\ell_1) f\\
 N(f)=\ell_1 e+ \ell_ 2 f\\
 N(g)=\frac{-1+\sqrt{1+4\a}}{2}p_2 e +p_2 f + (\frac{1+\sqrt{1+4\a}}{2\a}\ell_1+ k_1) g \\
 \end{array}
 \right.$.
 $(p_2\neq 0, \a > \frac{-1}{4}, \a \neq 0)$
  \end{center}\smallskip
 \end{enumerate}

 \item[(IIB)] $(p_1+p_2)^2-p_2(p_1+p_2)-\a p_2^2 \neq 0 $, then $\left\{
 \begin{array}{l}
 \ell_3 = k_3 =\ell_1= k_2= 0\\
 p_3= \ell_2=k_1\\
 p_2\neq 0\\
 (p_1+\frac{1}{2}p_2)^2-\frac{1+4\a}{4} p_2^2 \neq 0\\
 \end{array}
 \right.$.\smallskip
 \begin{enumerate}
 \item[(IIB1)] $\a<\frac{-1}{4}$, then $\left\{
 \begin{array}{l}
 \ell_3 = k_3 =\ell_1= k_2= 0\\
 p_3= \ell_2=k_1\\
 p_2\neq 0, \a<\frac{-1}{4}\\
 \end{array}
 \right.$.  That is to say,\smallskip
 \begin{center}$\left\{
 \begin{array}{l}
 N(e)=k_1 e\\
 N(f)=k_1 f\\
 N(g)=p_1 e +p_2 f + k_1 g \\
 \end{array}
 \right.~(p_2\neq 0, \a<\frac{-1}{4})$.
 \end{center}\smallskip
 \item[(IIB2)] $\a\ge \frac{-1}{4}$,
 %then $\left\{
 %\begin{array}{l}
 %\ell_3 = k_3 =\ell_1= k_2= 0\\
 %p_3= \ell_2=k_1\\
 %(p_1+ \frac{1+\sqrt{1+4\a}}{2}p_2)(p_1+ \frac{1-\sqrt{1+4\a}}{2}p_2)\neq 0\\
 %p_2\neq 0, \a\ge \frac{-1}{4}\\
 %\end{array}
 %\right.$.

 then $\left\{
 \begin{array}{l}
 \ell_3 = k_3 =\ell_1= k_2= 0\\
 p_3= \ell_2=k_1\\
 p_1\neq \frac{-1\pm \sqrt{1+4\a}}{2}p_2, p_2\neq 0, \a\ge \frac{-1}{4}, \a \neq 0\\
 \end{array}
  \right.$.
 That is to say,\smallskip
 \begin{center}$\left\{
 \begin{array}{l}
 N(e)=k_1 e\\
 N(f)=k_1 f\\
 N(g)=p_1 e +p_2 f + k_1 g \\
 \end{array}
 \right.~(p_1\neq \frac{-1\pm \sqrt{1+4\a}}{2}p_2, p_2\neq 0, \a\ge \frac{-1}{4}, \a \neq 0)$.
 \end{center}\smallskip
 \end{enumerate}
 \end{enumerate}
 \end{enumerate}

 \end{proof}

 \begin{thm}\label{thm:3nb12}  Let $(\ll, [,])$ be a Leibniz algebra and $[,]$ is given by
 \begin{center}
 \begin{tabular}{r|rrr}
 $[,]$ & $e$  & $f$ & $g$ \\
 \hline
 $e$ & $0$  & $0$ & $0$ \\
 $f$ & $0$  & $0$ & $0$  \\
 $g$ & $f$  & $0$ & $e$  \\
 \end{tabular}.
 \end{center}\vskip1mm
 Then all the Nijenhuis operators on $(\ll, [,])$ are given as follows:
 \begin{center}$~\left(
 \begin{array}{ccc}
 k_1 & k_2 & 0 \\
 0 & k_1 & 0 \\
 p_1 & p_2 & k_1 \\
 \end{array}
 \right).$
 \end{center} 
 \end{thm}

 \begin{proof} $N$ is Nijenhuis operator on $(\ll, [,])$ if and only if
 \begin{center}
 $\left\{
 \begin{array}{l}
 \ell_1 = \ell_3 = k_3= 0\\
 k_1 p_3+\ell_2^2-k_1 \ell_2-p_3 \ell_2=0\\
 p_1 p_3+ k_1 k_2+k_2 \ell_2-p_1 \ell_2-2 p_3 k_2=0 \\
 p_3 = k_1\\
 \end{array}
 \right.$.
 \end{center}
 Then $\left\{
 \begin{array}{l}
 \ell_1 = \ell_3 = k_3= 0\\
 p_3 = \ell_2 = k_1\\
 \end{array}
 \right.$. That is to say,\smallskip
 \begin{center}$\left\{
 \begin{array}{l}
 N(e)=k_1 e+k_2 f\\
 N(f)=k_1 f\\
 N(g)=p_1 e+ p_2 f+k_1 g\\
 \end{array}
 \right.$.
 \end{center}
 \end{proof}

 \begin{thm}\label{thm:3nb13}  Let $(\ll, [,])$ be a Leibniz algebra and $[,]$ is given by
 \begin{center}
 \begin{tabular}{r|rrr}
 $[,]$ & $e$  & $f$ & $g$ \\
 \hline
 $e$ & $0$  & $0$ & $0$ \\
 $f$ & $0$  & $0$ & $0$  \\
 $g$ & $e$  & $f$ & $0$  \\
 \end{tabular}.
 \end{center}\vskip1mm
 Then all the Nijenhuis operators on $(\ll, [,])$ are given as follows:
 \begin{center}$(1)~\left(
 \begin{array}{ccc}
 k_1 & k_2 & 0 \\
 \ell_1 & \ell_2 & 0 \\
 0 & 0 & p_3 \\
 \end{array}
 \right),$ \quad
 $(2)~\left(
 \begin{array}{ccc}
 k_1 & k_2 & 0 \\
 0 & \ell_2 & 0 \\
 0 & p_2 & \ell_2 \\
 \end{array}
 \right),~(p_2\neq 0),$  \quad
 \end{center}

 \begin{center}
 $(3)~\left(
 \begin{array}{ccc}
 k_1 & \frac{p_2^2 \ell_1}{p_1^2}+\frac{p_2 k_1}{p_1}-\frac{p_2 \ell_2}{p_1} & 0 \\
 \ell_1 & \ell_2 & 0 \\
 p_1 & p_2  & \frac{p_2 \ell_1}{ p_1} + k_1 \\
 \end{array}
 \right),~(p_1\neq 0).$
 \end{center} 
 \end{thm}

 \begin{proof} $N$ is Nijenhuis operator on $(\ll, [,])$ if and only if
 \begin{center}
 $\left\{
 \begin{array}{l}
 k_3 = \ell_3 = 0\\
 p_1(p_3-k_1)-p_2 \ell_1=0\\
 p_2(p_3-\ell_2)-p_1 k_2=0\\
 \end{array}
 \right.$.
 \end{center}
 \begin{enumerate}[(I)]
 \item $p_1=0$, then $\left\{
 \begin{array}{l}
 k_3 = \ell_3 = p_1 = 0\\
 p_2 \ell_1=0\\
 p_2(p_3-\ell_2)=0\\
 \end{array}
 \right.$. \smallskip
 \begin{enumerate}
 \item [(IA)] $p_2=0$, then $\left\{
 \begin{array}{l}
 k_3 = \ell_3 = p_1 = p_2 = 0\\
 \end{array}
 \right.$. That is to say,\smallskip
 \begin{center}$\left\{
 \begin{array}{l}
 N(e)=k_1 e + k_2 f\\
 N(f)=\ell_1 e + \ell_2 f\\
 N(g)=p_3 g\\
 \end{array}
 \right. $ .
 \end{center}\smallskip
 \item [(IB)] $p_2\neq 0$, we have
 $\left\{
 \begin{array}{l}
 k_3 = \ell_3 = p_1 = \ell_1 =0\\
 p_3 = \ell_2\\
 p_2\neq 0\\
 \end{array}
 \right. $.  That is to say,\smallskip
 \begin{center}$\left\{
 \begin{array}{l}
 N(e)=k_1 e+k_2 f\\
 N(f)=\ell_2 f\\
 N(g)=p_2 f+\ell_2 g\\
 \end{array}
 \right. ~(p_2\neq 0)$.
 \end{center}
 \end{enumerate}
 \item $p_1\neq 0$, then $\left\{
 \begin{array}{l}
 k_3 = \ell_3 = 0\\
 p_1\neq 0\\
 p_3= \frac{p_2 \ell_1}{ p_1} + k_1 \\
 k_2=\frac{p_2^2 \ell_1}{p_1^2}+\frac{p_2 k_1}{p_1}-\frac{p_2 \ell_2}{p_1} \\
 \end{array}
 \right.$. That is to say,\smallskip
 \begin{center}$\left\{
 \begin{array}{l}
 N(e)=k_1 e+(\frac{p_2^2 \ell_1}{p_1^2}+\frac{p_2 k_1}{p_1}-\frac{p_2 \ell_2}{p_1}) f\\
 N(f)=\ell_1 e+\ell_2 f\\
 N(g)=p_1 e+p_2 f+(\frac{p_2 \ell_1}{ p_1} + k_1) g\\
 \end{array}
 \right.~(p_1\neq 0)$.
 \end{center}
 \end{enumerate}
 \end{proof} 
 
 \section{Conclusion} We end this paper with one question. 
 %\begin{enumerate}[(I)] 
   %\item 
   In \cite{MSZ}, Ma, Sun and Zheng introduced the notion of Nijenhuis Leibniz bialgebras.
   \begin{defi}\label{de:yy} A {\bf Nijenhuis Leibniz bialgebra} is a vector space $\ll$ together with linear maps $[,]: \ll\o \ll\lr \ll$, $\delta : \ll\lr \ll\o \ll$, $N, S: \ll\lr \ll$ such that
   \begin{enumerate}[(1)]
    \item $(\ll,[,],\delta)$ is a Leibniz bialgebra.
    \item $(\ll,[,],N)$ is a Nijenhuis Leibniz algebra.
    \item $(\ll,\delta,S)$ is a Nijenhuis Leibniz coalgebra.\vskip-6mm
    \item For all $x, y\in \ll$, the equations below hold:
      \begin{eqnarray*}
       &S([N(x), y])+[x, S^2(y)]=[N(x), S(y)]+S([x, S(y)]),&\label{eq:c1}\\
       &S([x, N(y)])+[S^2(x), y]=[S(x), N(y)]+S([S(x), y]),&\label{eq:c2}\\
       &(S \o \id)\delta N+(\id \o N^2)\delta=(S\o N)\delta +(\id \o N)\delta N,&\label{eq:j1}\\
       &(\id \o S)\delta N+(N^2 \o \id)\delta=(N\o S)\delta +(N \o \id)\delta N.&\label{eq:j2}
     \end{eqnarray*}
   \end{enumerate}   
 Equation
 \begin{eqnarray*}\label{eq:cor171}
 &r_{12}r_{23}+r_{13}r_{23}=r_{12}^{\tau}r_{13}+r_{13}r_{12}^{\tau},&\label{eq:c171}
 \end{eqnarray*}
 where
 \begin{eqnarray*}
 r_{12}r_{23}=r^1\o [r^2, \bar{r}^1]\o \bar{r}^2,~~r_{13}r_{23}=r^1\o \bar{r}^1\o [r^2, \bar{r}^2],\\
 r_{12}r_{13}=[r^1, \bar{r}^1]\o r^2\o \bar{r}^2,~~r_{13}r_{12}=[r^1, \bar{r}^1]\o \bar{r}^2 \o r^2,
 \end{eqnarray*}
 and $\bar{r}=r,~~r^{\tau}=r^2\o r^1$, together with equations
 \begin{eqnarray*}
 &N(r^1)\o r^2=r^1\o S(r^2),&\label{eq:dks1}\\
 &S(r^1)\o r^2=r^1\o N(r^2)& \label{eq:dks2}
 \end{eqnarray*}
 is called an {\bf $S$-admissible classical Leibniz Yang-Baxter equation in $((\ll,[,]),  N)$} or simply an {\bf $S$-admissible cLYBe in $((\ll,[,]), N)$}.
 
 Let $(\ll, N)$ be an $S$-admissible Nijenhuis Leibniz algebra and $r:=r^1\o r^2\in \ll\o \ll$ be a symmetric solution of $S$-admissible cLYBe in $((\ll,[,]), N)$. Then $\big((\ll, N), \delta_r, S\big)$ is a Nijenhuis Leibniz bialgebra, where $\delta_r$ is given by
   $$
   \delta(x):=\delta_{r}(x)=-r^1\o [r^2, x]+[r^2, x]\o r^1+[x, r^2]\o r^1, \forall~x\in \ll.
   $$ In this case we call this bialgebra $\big((\ll, N), \delta_r, S\big)$ {\bf triangular}.
   \end{defi}
   \noindent In \cite{MSZ}, the authors also gave all the triangular Nijenhuis Leibniz bialgebra of dimension two. On the basis of all the Nijenhuis operators on the Leibniz algebras of dimension three given in Section \ref{se:3}, we can further discuss all the triangular Nijenhuis Leibniz bialgebra of dimension three. Or more generally, how can we derive the classifications of Nijenhuis Leibniz bialgebras of dimensions two and three?
   %\item
 %\end{enumerate}

 \end{document}